\baselineskip=17pt plus3pt minus1pt 
\magnification=\magstep1       
\hsize=5.5truein                     
\vsize=8.5truein                      
\parindent 0pt
\parskip=\smallskipamount
\mathsurround=1pt
\hoffset=.25truein
\voffset=.5truein
\nopagenumbers
%
%
\def\today{\ifcase\month\or
  January\or February\or March\or April\or May\or June\or
  July\or August\or September\or October\or November\or December\fi
  \space\number\day, \number\year}
%
 at 10truept
\font\smallheadfont=cmr8
\font\smallheadbf=cmbx8
\font\smallheadit=cmti8
%
\newcount\dispno      
\dispno=1\relax       
\newcount\refno       
\refno=1\relax        
\newcount\citations   
\citations=0\relax    
\newcount\sectno      
\sectno=0\relax       
\newbox\boxscratch    
%

%
%
%
\def\Section#1#2{\global\advance\sectno by 1\relax%
\label{Section\noexpand~\the\sectno}{#2}%
\smallskip
\goodbreak
\setbox\boxscratch=\hbox{\bf Section \the\sectno.~}%
{\hangindent=\wd\boxscratch\hangafter=1
\noindent{\bf Section \the\sectno.~#1}\nobreak\smallskip\nobreak}}
%
\def\sqr#1#2{{\vcenter{\vbox{\hrule height.#2pt
              \hbox{\vrule width.#2pt height#1pt \kern#1pt
              \vrule width.#2pt}
              \hrule height.#2pt}}}}
\def\square{$\mathchoice\sqr34\sqr34\sqr{2.1}3\sqr{1.5}3$}
\def\endproof{~~\hfill\square\par\medbreak}
\def\noproof{~~\hfill\square}
%
%
\def\proc#1#2#3{{\hbox{${#3 \subseteq} \kern -#1cm _{#2 /}\hskip 0.05cm $}}}
\def\propcont{\mathchoice\proc{0.17}{\scriptscriptstyle}{}
                         \proc{0.17}{\scriptscriptstyle}{}
                         \proc{0.15}{\scriptscriptstyle}{\scriptstyle }
                         \proc{0.13}{\scriptscriptstyle}{\scriptscriptstyle}}
%

%
\def\normalin{\hbox{\raise0.045cm \hbox
                   {$\underline{\triangleleft }$}\hskip0.02cm}}
%
%
\def\'#1{\ifx#1i{\accent"13 \i}\else{\accent"13 #1}\fi}
%
%
%
\def\semidirect{\rlap{$\times$}\kern+7.2778pt \vrule height4.96333pt
width.5pt depth0pt\relax\;}
%
%
\def\prop#1#2{\noindent{\bf Proposition~\the\sectno.\the\dispno. }%
\label{Proposition\noexpand~\the\sectno.\the\dispno}{#1}\global\advance\dispno 
by 1{\it #2}\smallbreak}
\def\thm#1#2{\noindent{\bf Theorem~\the\sectno.\the\dispno. }%
\label{Theorem\noexpand~\the\sectno.\the\dispno}{#1}\global\advance\dispno
by 1{\it #2}\smallbreak}
\def\cor#1#2{\noindent{\bf Corollary~\the\sectno.\the\dispno. }%
\label{Corollary\noexpand~\the\sectno.\the\dispno}{#1}\global\advance\dispno by
1{\it #2}\smallbreak}
\def\defn{\noindent{\bf
Definition~\the\sectno.\the\dispno. }\global\advance\dispno by 1\relax}
\def\lemma#1#2{\noindent{\bf Lemma~\the\sectno.\the\dispno. }%
\label{Lemma\noexpand~\the\sectno.\the\dispno}{#1}\global\advance\dispno by
1{\it #2}\smallbreak}
\def\rmrk#1{\noindent{\bf Remark~\the\sectno.\the\dispno.}%
\label{Remark\noexpand~\the\sectno.\the\dispno}{#1}\global\advance\dispno
by 1\relax}
\def\proof{\noindent{\it Proof: }}
\def\numbeq#1{\the\sectno.\the\dispno\label{\the\sectno.\the\dispno}{#1}%
\global\advance\dispno by 1\relax}

\def\comm#1,#2{\left[#1{,}#2\right]}
\newdimen\boxitsep \boxitsep=0 true pt
\newdimen\boxith \boxith=.4 true pt 
\newdimen\boxitv \boxitv=.4 true pt
\gdef\boxit#1{\vbox{\hrule height\boxith
                    \hbox{\vrule width\boxitv\kern\boxitsep
                          \vbox{\kern\boxitsep#1\kern\boxitsep}%
                          \kern\boxitsep\vrule width\boxitv}
                    \hrule height\boxith}}
\def\square{\ \hbox{\vrule height7.5pt depth1.5pt width 6pt}\par}
\outer\def\square{\ifmmode\else\hfill\fi
   \setbox0=\hbox{} \wd0=6pt \ht0=7.5pt \dp0=1.5pt
   \raise-1.5pt\hbox{\boxit{\box0}\par}
}

\def\frac#1/#2{\leavevmode\kern.1em
              \raise.5ex\hbox{\the\scriptfont0 #1}\kern-.1em
              /\kern\.15em\lower.25ex\hbox{\the\scriptfont0 #2}}
\def\incnoteq{\lower.1ex \hbox{\rlap{\raise 1ex
     \hbox{$\scriptscriptstyle\subset$}}{$\scriptscriptstyle\not=$}}}
%
%
\def\mapright#1{\smash{
     \mathop{\longrightarrow}\limits^{#1}}}
\def\mapdown#1{\Big\downarrow
 \rlap{$\vcenter{\hbox{$\scriptstyle#1$}}$}}


\def\propcontup{\bigcup\!\!\!\rlap{\kern+.2pt$\backslash$}\,\kern+1pt\vert}
%
%
%
\def\label#1#2{\immediate\write\aux%
{\noexpand\def\expandafter\noexpand\csname#2\endcsname{#1}}}
%
\def\ifundefined#1{\expandafter\ifx\csname#1\endcsname\relax}
%
%
\def\ref#1{%
\ifundefined{#1}\message{! No ref. to #1;}%
 \else\csname #1\endcsname\fi}
%
%
\def\refer#1{%
\the\refno\label{\the\refno}{#1}%
\global\advance\refno by 1\relax}
%
%
\def\cite#1{%
\expandafter\gdef\csname x#1\endcsname{1}%
\global\advance\citations by 1\relax
\ifundefined{#1}\message{! No ref. to #1;}%
\else\csname #1\endcsname\fi}
%
%
\font\bb=msbm10 
\font\bp=msbm10 at 8truept      
%
%
%

\def\Q{\hbox{\bb Q}}

\def\Z{\hbox{\bb Z}}                     \def\ZZ{\hbox{\bp Z}}

\def\Z{\hbox{\bb Z}}                     \def\ZZ{\hbox{\bp Z}}

\baselineskip=12pt
\parskip=4pt plus 1pt minus 1pt
\newread\aux
\immediate\openin\aux=\jobname.aux
\ifeof\aux \message{! No file \jobname.aux;}
\else \input \jobname.aux \immediate\closein\aux \fi
\newwrite\aux
\immediate\openout\aux=\jobname.aux

\font\smallheadfont=cmr8
\font\smallheadbf=cmbx8
\font\smallheadit=cmti8

\null\vskip 1truein
\centerline{DOMINIONS IN VARIETIES OF NILPOTENT GROUPS}
\vskip \baselineskip
\vskip \baselineskip
\centerline{Arturo Magidin}
\vskip\baselineskip
{\baselineskip 12pt plus2pt minus1pt
\leftskip 2truein Instituto de Matem\'aticas, Oficina 112\par
\leftskip 2truein UNAM, Circuito Exterior, Cd. Universitaria\par
\leftskip 2truein 04510 M\'exico City, MEXICO\par
\leftskip 2truein e-mail: magidin@matem.unam.mx\par}
\vskip \baselineskip
\vskip \baselineskip
\vskip \baselineskip

{\parindent=20pt
\narrower\narrower

\noindent{Abstract. We investigate the concept of
dominion (in the sense of Isbell) in several varieties
of nilpotent groups. We obtain a complete description of dominions in
the variety of nilpotent groups of class at most~two.
Then we look at the behavior of dominions of
subgroups of groups in ${\cal N}_2$ when taken in the
context of ${\cal N}_c$ for $c>2$.
Finally, we establish the existence of nontrivial
dominions in the category of all nilpotent groups.}}

\bigskip
\medskip



\Section{Introduction}{intro}

Suppose a group~$G$ and a subgroup~$H$ of~$G$ are given. Are there any
elements $g\in G\setminus H$ such that any two group morphisms
which agree on~$H$ must also agree on~$g$? What if we require all the
groups involved to be~nilpotent?

To put this question in a more general context, let~$\cal C$ be a full
subcategory of the category of all algebras (in the sense of Universal
Algebra) of a fixed type, which is closed under passing to subalgebras.
Let $A\in {\cal C}$, and let~$B$ be a subalgebra of~$A$. Recall that,
in this situation, Isbell~{\bf [\cite{isbellone}]} defines the {\it
dominion of~$B$ in~$A$} (in the category ${\cal C}$) to be the
intersection of all equalizer subalgebras of~$A$
containing~$B$. Explicitly,
$${\rm dom}_A^{\cal C}(B)=\Bigl\{a\in A\bigm| \forall C\in {\cal C}\;
\forall f,g\colon A\to C,\ {\rm if}\ f|_B=g|_B{\rm\ then\ }
f(a)=g(a)\Bigr\}.$$

Therefore, the question with which we opened this discussion may be
restated in terms of the dominion of~$H$ in~$G$ in the category of~context.

If ${\rm dom}_A^{\cal C}(B)=B$ we will say that the dominion of~$B$
in~$A$ is {\it trivial\/}, and we will say it is {\it
nontrivial\/}~otherwise.

In the case where the category~${\cal C}$ is actually a {\sl variety}
(or more generally, a right closed category; see {\bf
[\cite{isbellone}]}), we can look at the amalgamated coproduct of two
copies of~$A$, amalgamated over the subalgebra~$B$. This is the
pushout of the diagram {$$
\matrix{B&\mapright{i}&A\cr
\mapdown{i}\cr
A\cr}$$}\relax
and is denoted by $A\amalg_B^{\cal C} A$. If we write $(\lambda,\rho)$
for the universal pair of maps from~$A$ to $A\amalg_B^{\cal C} A$, then ${\rm
dom}_A^{\cal C}(B)$ is the equalizer of~$\lambda$ and~$\rho$; one can
in fact verify that ${\rm dom}_A^{\cal
C}(B)=\lambda(A)\cap\rho(A)$. By a classical theorem of Schreier, if
${\cal C}={\cal G}roup$, then for all $A\in {\cal C}$ and for every
subgroup $B$ of~$A$, ${\rm dom}_A^{{\cal C}}(B)=B$.

Dominions are related to group amalgams, and particularly to {\it
special amalgams}. Recall that an amalgam of $A$ and~$C$ with core~$B$
consists of groups $A$, $C$, and~$B$, equipped with one-to-one morphisms
$\Phi_A\colon B\to A$, $\Phi_C\colon B\to C$. We will denote this
situation by writing $[A,C;B]$. We say that the
amalgam is {\it weakly embeddable} in ${\cal C}$ if there exists a
group~$M\in {\cal C}$ and one-to-one mappings
$$\lambda_A\colon A\to M,\qquad \lambda_C\colon C\to M,\qquad
\lambda_B\colon B\to M$$
such that
$$\lambda_A\circ\Phi_A=\lambda_B,\qquad\qquad
\lambda_C\circ\Phi_C=\lambda_B.$$
We usually identify~$B$ with its images in~$A$ and~$C$.

We say the amalgam is {\it strongly embeddable} if, furthermore, there
is no identification between elements of $A\setminus B$ and
$C\setminus B$.
Finally, by a {\it special amalgam} we mean an amalgam $[A,C;B]$,
where there is an isomorphism $\alpha\colon A\to C$ such that
\hbox{$\alpha\circ\Phi_A=\Phi_C$}. In this case, we usually write $[A,A;B]$,
with $\alpha={\rm id}_A$ being~understood.

It is not hard to see that $B$ equals its own dominion in~$A$ (in the
category~${\cal C}$) if and only if the special amalgam $[A,A;B]$ is
strongly embeddable, and that the dominion of~$B$ is in general the
least subgroup $D$ of~$A$, such that~$D$ contains~$B$, and $[A,A;D]$
is strongly embeddable. 
We refer the reader to the survey article by
Higgins {\bf [\cite{episandamalgs}]} for the~details.

Dominions are also related to epimorphisms, and in fact were
introduced by Isbell to study them. Recall that given a category of
algebras~${\cal C}$, a map $f\colon A\to C$ is an epimorphism if and
only if it is right cancellable in~${\cal C}$. That is, if for all
pairs of maps $g,h\colon C\to K$ in~${\cal C}$, $g\circ f=h\circ f$
implies that $g=h$.  Clearly, $f\colon A\to C$ is
an epimorphism in~${\cal C}$ if and only if ${\rm dom}_C^{\cal
C}(f(A))=C$.

We note that in categories of algebras, all surjective maps are
epimorphisms; in some categories, such as the category of all groups,
the converse also holds. But the converse does not hold in general:
for example, the embedding $\Z\hookrightarrow\Q$ is an epimorphism in
the category of all~rings.

Unfortunately, the connection to epimorphisms is not relevant in the
context of the present work. Peter Neumann has proven {\bf
[\cite{pneumann}]} that in all ``reasonable'' categories of solvable
groups, all epimorphisms are necessarily surjective; specifically, he
showed that in any full category which consists of solvable groups and is
closed under taking quotients, all epimorphisms are surjective. This
was later extended substantially by S. McKay {\bf [\cite{mckay}]}.

Despite these negative results, one should not set aside dominions as
worthless in the context of varieties of nilpotent groups. Their
relation to amalgams can be useful; for example, the work on dominions
in ${\cal N}_2$ done here is used elsewhere {\bf [\cite{absclosed}]}
to characterize the special amalgamation bases in that variety. Also,
from a very general point of view, finding nontrivial dominions should
be enough to be of interest. Any result that says that we can predict
the behavior of a function at some point based on partial information
(in this case, the value of the morphism at points in the dominion
of~$H$ but not in~$H$, based on the value on~$H$) has the potential of
being~useful. Aside from this, the dominion construction determines a
special class of subalgebras, those which equal their own dominion;
equivalently, those which are ``closed'' under the closure operation
induced by the construction. When not all subalgebras are closed, that
is when not every subalgebra equals its own dominion, it is possible
for this class to have interesting properties of their own. For
example, Bergman has shown {\bf [\cite{ordersberg}]} that in the
category of Orderable Groups (groups in which an order can be defined
which is compatible with the operations), the class of dominion-closed
subgroups of an orderable group~$G$ are precisely the subgroups~$T$
for
which the amalgamated coproduct $G\amalg_T^{{\cal G}roup} G$ is
also~orderable.
 
In \ref{basicres} we prove some of the basic properties of dominions, and
recall the basic definitions associated to nilpotent groups. In
\ref{domsniltwo} we will study the variety of nilpotent groups of
class at most two. Then, in \ref{domstwoengel} we will generalize the
arguments in \ref{domsniltwo} to the variety of 2-Engel groups. 
In \ref{domsniltwoinmetab} we will study how dominions of subgroups of
${\cal N}_2$-groups behave when we change the variety of context
from~${\cal N}_2$ to ${\cal A}^2 \cap {\cal N}_c$, for $c>2$; and in
\ref{domsmetabnil} we will look at the category of all metabelian
nilpotent groups. In that section, we will have some things to say
about dominions in the variety of all metabelian groups as well. Then, in
\ref{domsniltwoinnilk} and \ref{nontrivnil} we will expand these
investigations to cover the case where the category of context is
${\cal N}_c$, and ${\cal N}\!il$, respectively.
Finally, in \ref{domsinothernil} we will mention some related results
in other varieties of nilpotent~groups.

The contents of this work are part of the author's doctoral
dissertation, which was conducted under the direction of Prof.~George
M.~Bergman, at the University of California at~Berkeley. It is my very
great pleasure to express my deep gratitude and indebtedness to
Prof.~Bergman for his advice and encouragement; his many suggestions
have improved this work in ways too numerous to list explicitly; he
also caught and helped correct many mistakes. Any errors that remain,
however, are entirely my own responsibility. 

\bigskip

\Section{Preliminary definitions and results}{basicres}

The group operation will be written multiplicatively unless otherwise
stated; given a group~$G$, the identity element of~$G$ will be denoted
by~$e_G$, with the subscript omitted if it is understood from
context. Given two elements $x$ and~$y$ in~$G$, we write $x^y=y^{-1}xy$, and
we will denote their commutator by $[x,y]=x^{-1}y^{-1}xy$. Given two
subsets $A,B$ of~$G$ (not necessarily subgroups), we denote by $[A,B]$
the subgroup of~$G$ generated by all elements $[a,b]$ with $a\in A$
and $b\in B$.  We also define inductively the left-normed commutators
of weight~\hbox{$c+1$}:
$$[x_1,\ldots,x_c,x_{c+1}] =
\bigl[[x_1,\ldots,x_c],x_{c+1}\bigr];\quad c\geq 2.$$

The centralizer in~$G$ of a subgroup~$H$ will be denoted by
$C_G(H)$. We will denote the center of~$G$ by $Z(G)$. 

A variety of groups is a full subcategory of ${\cal G}roup$ which is
closed under taking subgroups, quotients, and arbitrary direct products.
We will first establish some basic properties of dominions in
varieties (and more general categories) of~groups.

\lemma{easyone}{Let ${\cal C}$ be a category of groups.}
{\parindent=20pt\it
\item{(i)} If $A\in{\cal C}$, then ${\rm dom}_A^{\cal C}({-})$ is a
closure operator on the lattice of subgroups of~$A$.
\item{(ii)} Given a homomorphism $h\colon A\to A'$, with $A,A'\in
{\cal C}$, if~$B$ is a subgroup of~$A$, then
$h\bigl({\rm dom}_A^{\cal C}(B)\bigr)\subseteq {\rm dom}_{A'}^{\cal
C}\bigl(h(B)\bigr)$.
\item{(iii)} If $A,C\in {\cal C}$ and $B,C$ are subgroups of~$A$,
with $B\subseteq C$, then
$${\rm dom}_{C}^{\cal C}(B)\subseteq {\rm dom}_A^{\cal C}(B).$$
\par}
 
\proof (i) and (ii) are immediate from the definition of
dominion; (iii) follows from~(ii) by considering the inclusion map
$i\colon C\to A$.\endproof

\lemma{varinclusions}{Let ${\cal C}\subseteq{\cal D}$ be two
full subcategories of groups, and let $G\in{\cal C}$. If $H$ is a
subgroup of~$G$, then
${\rm dom}^{\cal D}_{G}(H) \subseteq {\rm dom}^{\cal C}_{G}(H).$}
 
\proof  The set of pairs of maps from~$G$ to groups
in~${\cal C}$ is a subset of the set of pairs of maps from~$G$ to
groups in~${\cal D}$. The inclusion now follows from the definition of
the dominion as the intersection of equalizer subgroups
of~$G$.\endproof
 
\lemma{easysix}{Let ${\cal C}$ be a full subcategory of~${\cal
G}roup$, and let $G\in {\cal C}$ be a group. If $H$ and $K$ are
subgroups of~$G$, then $$\Bigl\langle{\rm dom}_{G}^{\cal C}({H}), {\rm
dom}_{G}^{\cal C}({K})\Bigr\rangle
\subseteq {\rm dom}_{G}^{\cal C}\Bigl(\bigl\langle
H,K\bigr\rangle\Bigr).$$}
 
\proof Let $L=\langle H,K\rangle$. Recall that ${\rm dom}_G^{\cal
C}(-)$ is a closure operator on the
subgroups of~$G$. Now, $H\subseteq L$ implies that the dominion of~$H$
is contained in the dominion of~$L$; the same is true of the dominion
of~$K$, and therefore it is also true of the subgroup that the
dominion of~$H$ and the dominion of~$K$ generate, establishing
the~result.\endproof
 
\thm{dirproduct}{Let $G_1$ and $G_2$ be groups, and let $H_1$ be a
subgroup of $G_1$, $H_2$ a subgroup of $G_2$. If ${\cal C}$ is any
full subcategory of~${\cal G}roup$ such
that $G_1$, $G_2$, and~$G_1\times G_2$ are in~${\cal C}$, then
$${\rm dom}_{G_1\times G_2}^{\cal C}({H_1\times H_2}) = {\rm
dom}_{G_1}^{\cal C}({H_1})\times
{\rm dom}_{G_2}^{\cal C}({H_2}).$$}
 
\proof We can identify $G_1$ with the subgroup $G_1\times \{e\}$ of
$G_1\times G_2$, and similarily $G_2$ with the subgroup $\{e\}\times
G_2$. It follows, from
\ref{easysix}, that ${\rm dom}_{G_1}^{\cal C}({H_1})$ and ${\rm
dom}_{G_2}^{\cal C}({H_2})$ are
contained in ${\rm dom}_{G_1\times G_2}^{\cal C}({H_1\times
H_2})$. Therefore
$$\eqalign{{\rm dom}_{G_1\times G_2}^{\cal C}({H_1\times
H_2})&\supseteq \Bigl\langle
{\rm dom}_{G_1}^{\cal C}({H_1}),{\rm dom}_{G_2}^{\cal
C}({H_2})\Bigr\rangle\cr
&= {\rm dom}_{G_1}^{\cal C}({H_1})\times{\rm dom}_{G_2}^{\cal
C}({H_2}).\cr}$$
Now let $(g_1,g_2)\notin
{\rm dom}_{G_1}^{\cal C}({H_1})\times{\rm dom}_{G_2}^{\cal
C}({H_2})$. Assume, without loss of
generality, that $g_1\notin {\rm dom}_{G_1}^{\cal C}({H_1})$.
 
Therefore there is a group $K\in{\cal C}$, and a pair of maps
$\psi,\phi\colon G_1\to K$ such that $\psi|_{H_1}=\phi|_{H_1}$, but
$\psi(g_1)\not= \phi(g_1)$. Let $\pi\colon G_1\times G_2\to G_1$ be
the canonical projection, and compare the maps $\psi\circ\pi$ with
$\phi\circ\pi$. By construction, they agree on $H_1\times H_2$, but
disagree on $(g_1,g_2)$. Therefore, $(g_1,g_2)\notin {\rm
dom}_{G_1\times
G_2}^{\cal C}({H_1\times H_2})$.\endproof
 
\rmrk{notinfinite} We add a caution, however. \ref{dirproduct} implies
that the same result holds for a finite number of direct factors, and
that the analogous result holds for the direct sum of an arbitrary
number of factors. However, in the case of an infinite direct product,
equality may no longer~hold. An example of this will be given below,
in Example~\ref{counterexone}.
 
\lemma{easyfour}{If ${\cal C}$ is a full subcategory of groups which
is closed under quotients, then normal subgroups are
dominion-closed. That is, if~$N$ is a
normal subgroup of~$G$, with $G\in {\cal C}$, then ${\rm
dom}_{G}^{{\cal C}}(N)= N$.}
 
\proof Compare the maps $\pi,\zeta\colon G\to G/N$, where $\pi$ is
the canonical epimorphism onto the quotient and $\zeta$ is the zero
map. They both agree on $N$, and disagree on any element not
in~$N$.\endproof
 
\cor{normalclosure}{Let ${\cal C}$ be a full subcategory of groups
which is closed under quotients. If $G$ is a group in ${\cal C}$
and~$H$ is a
subgroup of~$G$, then
$${\rm dom}_G^{\cal C}(H)\subseteq H[H,G].$$}
 
\proof The subgroup $H[H,G]$ is the normal closure of the
subgroup~$H$; hence it contains~$H$, and is normal in~$G$. By
\ref{easyfour}, 
$${\rm dom}_G^{\cal C}(H)\subseteq {\rm dom}_G^{\cal
C}(H[H,G])=H[H,G].\eqno\noproof$$ 
 
\rmrk{abvar} It follows from \ref{easyfour} that dominions in any
variety ${\cal V}$ of abelian groups are trivial (that is,
${\rm
dom}_{G}^{{\cal V}}({H})=H$ for all groups $G$ and subgroups~$H$ in
the variety), since any subgroup is~normal. In fact, we see that for
any variety $\cal V$, if $G\in \cal V$ is abelian, then ${\rm
dom}^{\cal V}_{G}(H)=H$ for any subgroup~$H$ of~$G$. We express this
situation by saying that {\it dominions of subgroups of abelian groups
are trivial (in any variety ${\cal V}$).\/}

Recall that a class of groups~${\cal P}$ is a {\it pseudovariety} if
it is closed under quotients, subgroups, and {\sl finite}
direct~products.
 
\prop{quotient}{Let $\{{\cal P}_i\}_{i\in I}$ be a nonempty
collection of pseudovarieties. Let $G\in \bigcup{\cal P}_i$, $H$ a
subgroup of~$G$, and $N\triangleleft G$ such that $N\subseteq H$. Then
$${\rm dom}_{G/N}^{{\cup}{\cal P}_i}(H/N) = {\rm dom}_G^{{\cup}{\cal
P}_i}(H) \Bigm/ N.$$}
 
\proof  We have that ${\rm dom}_G^{{\cup}{\cal P}_i}(H) \,/\,N
\subseteq {\rm dom}_{G/N}^{{\cup}{\cal P}_i}(H/N)$, using
\ref{easyone}(ii) and setting $h\colon
G\to G/N$ to be the canonical surjection.
 
For the reverse inclusion, assume $x\notin {\rm dom}_G^{{\cup}{\cal
P}_i}(H)$. Then there exists a group $K\in \bigcup {\cal P}_i$
and a pair of maps $f,g\colon G\to K$ such that $f|_H=g|_H$ and
$f(x)\not=g(x)$.
 
Consider the induced homomorphisms $(f\times f), (f\times g)\colon
G\to K\times
K$, and let $L$ be the subgroup of $K\times K$ generated by the images
of~$G$ under these two morphisms. Since~$N$ is normal in~$G$, and
contained in~$H$, the common image of~$N$ under these two maps
will be normal in~$L$. This image may be written as the set
$$(f\times f)(N)=\Bigl\{ \bigl(f(n),f(n)\bigr) \in K\times K\,\Bigm|
\,
n\in N\Bigr\}.$$
 
We claim that $(f\times f)(x)$ and $(f\times g)(x)$ are not in the
same coset of
$(f\times f)(N)$ in~$L$. This will prove the claim, since we can then
mod out
by~$(f\times f)(N)$ to obtain an induced pair of maps $G/N\to
L/(f\times f)(N)$ which agree on $H/N$ and disagree on~$xN$.
 
Indeed, $\bigl((f\times f)(x)\bigr)\bigl((f\times g)(x)\bigr)^{-1} =
(e,f(x)g(x)^{-1})$, and we know that
the second coordinate is not trivial, because $f$ and~$g$ disagree
on~$x$. Hence,
$$\bigl((f\times f)(x)\bigr) \bigl((f\times
g)(x)\bigr)^{-1}$$ is not a diagonal element,
and in particular cannot lie in $(f\times f)(N)$. This proves
the~proposition.\endproof
 
\rmrk{whatitmeans} Note that a collection of groups is a union of
pseudovarieties if and only if it contains the pseudovariety generated
by each of its members. In particular, a collection of groups is a
union of pseudovarieties if and only if it is closed under subgroups,
quotients, and {\it squares}, where the square of a group~$G$
is the group $G\times G$.
 
\prop{centralizer}{{\rm (Cf. Corollary~2.4 in~{\bf
[\cite{isbellone}]})} Let ${\cal C}$ be a full subcategory of~${\cal
G}roup$, and let $G\in {\cal C}$. If $H$
is a subgroup of~$G$, then
$$C_G(H) = C_G({\rm dom}_{G}^{\cal C}({H})).$$}
 
\proof Since $H\subseteq{\rm dom}_{G}^{\cal C}({H})$, we automatically
have
$$C_G({\rm dom}_{G}^{\cal C}({H}))\subseteq C_G(H).$$
 
To prove the reverse inclusion, let $g\in C_G(H)$ and consider the
inner automorphism $\phi_g$ of~$G$ given by conjugation
by~$g$. Since $g\in C_G(H)$, it follows that $\phi_g$ fixes $H$
pointwise. Therefore $\phi_g|_H = {\rm id}_G|_H$, so it follows
that
$\phi_g$ also fixes ${\rm dom}_{G}^{\cal C}({H})$. That is, for all
$d$ in
${\rm dom}_{G}^{\cal C}({H})$,
$d=\phi_g(d)=g^{-1}dg$, hence $gd=dg$. So $g$ lies in $C_G({\rm
dom}_{G}^{\cal C}({H}))$,
as~claimed.\endproof

\cor{domofab}{{\rm (Cf. Cor 2.5 in~{\bf [\cite{isbellone}]})} Let ${\cal C}$
be a full subcategory of ${\cal G}roup$,
and let $G$ be a group in~${\cal C}$. If~$H$ is an
abelian
subgroup of~$G$, then ${\rm dom}_{G}^{\cal C}({H})$ is also~abelian.}
 
\proof By \ref{centralizer}, since all elements of~$H$
centralize~$H$, they also centralize ${\rm dom}_G^{\cal
C}(H)$. Therefore, every element of the dominion of~$H$ commutes with
every element of~$H$; hence every element of the dominion
centralizes~$H$, and therefore also centralizes the~dominion
of~$H$.\endproof

Next we recall the basic definitions and terminology associated to
nilpotent groups. We refer the reader to~{\bf
[\cite{hneumann}]} and~{\bf [\cite{rotman}]} for the~proofs.

\defn For a group~$G$ we define the {\it lower central series} of~$G$
recursively as follows: $G_1 = G$, and $G_{c+1} = [G_c,G]$ for $c\geq
1$. We call $G_c$ the $c$-th term of the lower central series of~$G$; $G_c$
is generated by elements of the form $[x_1,\ldots,x_c]$ for $c\geq
2$, and $x_i$ ranging over the elements of~$G$. 
We sometimes also write $G'=G_2=[G,G]$ for the commutator subgroup of~$G$.

The factor groups $G_{i-1}/G_i$ are called the
{\it lower central factors} of~$G$.

A group~$G$ is {\it nilpotent of class $c$} if and only if
$G_{c+1}=\{e\}$.

We will denote by ${\cal A}$ the variety of all abelian groups; by
${\cal A}^2$ the variety of all metabelian groups (that is, solvable
groups of solvability length at most~two);
by ${\cal N}_c$ the variety
of all nilpotent groups of class at most~$c$; by ${\cal A}^2\cap
{\cal N}_c$ the variety of all metabelian nilpotent groups of class at
most~$c$; by ${\cal B}_n$ the Burnside variety of exponent~$n$,
consisting of all groups that satisfy the identity~$x^n=e$;
by ${\cal N}\!il$ the category of all nilpotent groups,
and by ${\cal A}^2\cap {\cal N}\!il$ the category of all
met\discretionary{a-}{a}{a}belian nilpotent groups. Note that the last
two classes are not varieties, since they are not closed under
arbitrary direct~products.

The following lemma, which is easily established by direct
computation, and the definitions that follow it, will be useful in
subsequent considerations.

\lemma{commident}{The
following hold for any elements $x$,~$y$, $z$, and~$w$  of
an arbitrary group~$G$:}
{\parindent=30pt\it
\item{(a)}{$xy = yx[x,y];\qquad x^y=x[x,y]$.}
\item{(b)}{$[x,y]^{-1} = [y,x]$.}
\item{(c)}{$[xy,z] = [x,z]^y[y,z] = [x,z][x,z,y][y,z]$.}
\item{(d)}{$[x,zw] = [x,w][x,z]^w = [x,w][x,z][x,z,w]$.\noproof}\par}

\defn Let $G$ be a group, generated by elements
$x_1,\ldots,x_n$. We define the {\it basic commutators} and the
ordering among them recursively, as~follows:\par
{\parindent=30pt
\item{(i)}{The letters $x_1,\ldots,x_n$ are basic commutators of
weight one, ordered by setting $x_i<x_j$ if $i<j$.}\par
\item{(ii)}{If the basic commutators $c$ of weight less than~$k$
have been 
defined and ordered, define the basic commutators of weight $k$ by the
rules: $[c_1,c_2]$ is a basic commutator of weight $k$ if}\par
\itemitem{(a)}{${\rm weight}(c_1) + {\rm weight}(c_2) = k$,}\par
\itemitem{(b)}{$c_1>c_2$,}\par
\itemitem{(c)}{If $c_1 = [c_3,c_4]$, then $c_2\geq c_4$.}\par
\item{(iii)} Then continue the order by setting $c>c'$ if ${\rm
weight}(c)>{\rm weight}(c')$, and ordering those of order $k$
lexicographically; explicitly, we set $[c_1,c_2]<[c'_1,c'_2]$ if
$c_1<c'_1$, or if $c_1=c'_1$ and $c_2<c'_2$.\par}
 
\thm{generalcollect}{{\rm (Propositions~31.52 and~31.53 in~{\bf
[\cite{hneumann}]})} If 
$G=\langle g_1,\ldots,g_n\rangle$ is nilpotent of class~$c$, and
$\alpha\colon F(x_1,\ldots,x_n)\to G$ is the map from the relatively
free nilpotent group of class $c$ and rank $n$ sending $x_i$ to $g_i$,
then every element $g\in G$, $g\not= e$ is a product $g=\alpha\bigl(
c_{1}^{m_1}\cdots c_{\ell}^{m_\ell}\bigr)$, where the $c_{\lambda}$
are basic commutators of weight less than or equal to $c$,
$c_{1}<\cdots<c_{\ell}$, and the $m_{\lambda}$ are
nonzero~integers. If $G$ is freely generated in ${\cal N}_c$ by the
$g_i$, then the $m_{\lambda}$ take independently all nonzero integral
values, and the representation is~unique.\noproof}
 
\defn All elements of a group~$G$ are said to be {\it
commutators of weight~$1$}. The {\it commutators of weight~$n$} are
defined recursively as the elements $[x,y]$ such that $x$ is a
commutator of weight~$k$, $y$ is a commutator of weight~$m$, and
$k+m=n$.

Note that an element $g\in G$ may be considered to be a commutator of
several different weights.

\lemma{commweights}{{\rm (Lemma~33.35 in~{\bf [\cite{hneumann}]})} Every
commutator of weight $n$ of a group $G$ is a
product of left-normed commutators of weight~$n$ and their inverses, and it
belongs to $G_{n}$, the $n$-th term of the lower central series
of~$G$.\noproof}

\goodbreak
\Section{Dominions in ${\cal N}_2$}{domsniltwo}

A group~$G$ is nilpotent of class at most two if and only if
$G'\subseteq Z(G)$; that is, if and only if commutators
are~central. From this, it is easy to show the following two results:

\prop{commidenttwo}{Let $G\in {\cal N}_2$. Then $\forall x,y,z \in G$
$$[xy,z] = [x,z][y,z],\qquad [x,yz] =
[x,y][x,z].\eqno{\noproof}$$}

\cor{commpowerform}{Let $G\in {\cal N}_2$. Then for all $x,y\in G$ and
for all $n\in\Z$, \hbox{$[x^n,y]=[x,y]^n=[x,y^n]$}.\noproof}

It is not hard to verify that the group~$F(x,y)$ presented by
$$F(x,y)=\langle x,y\,|\, [x,y,x]=[x,y,y]=e\rangle.$$
is the relatively free ${\cal N}_2$-group of rank two, freely
generated by~$x$ and~$y$.

Let $F=F(x_1,\ldots,x_n)$ be the relatively free ${\cal N}_2$-group in
$n$ generators. In~$F$ there is a normal form for the
elements. Namely, every element can be written uniquely in the form
$x_1^{m_1}x_2^{m_2}\cdots x_n^{m_n}u$
where $u$ is of the form
$$u = \prod_{1\leq i < j \leq m}[x_j,x_i]^{m_{ji}}\qquad
m_{ij}\in\hbox{\Z};$$
see~{\bf [\cite{metab}]} for the~details.

Furthermore, if $g=x_1^{m_1}\cdots x_n^{m_n}u$ and $g'=x_1^{p_1}\cdots
x_n^{p_n}u'$, then $$gg'=x_1^{m_1+p_1}\cdots
x_n^{m_n+p_n}v\leqno(\numbeq{prodformulainntwo})$$ where
$$v=uu'\prod_{1\leq i < j \leq
n}[x_j^{m_j},x_i^{p_i}]
             =uu'\prod_{1\leq i < j \leq
               n}[x_j,x_i]^{m_jp_i}$$
and
$$g^{-1} = x_1^{-m_1}\cdots x_n^{-m_n}u^{-1}\prod_{1\leq i<j\leq n}
[x_j,x_i]^{m_jm_i}.\leqno(\numbeq{inverseformulaintwo})$$

\lemma{description}{Let $F(x_1,\ldots,x_n)$ be the relatively free
group on $n$ generators in~${\cal N}_2$. Let $H=\langle
x_1^{a_1},\ldots,x_n^{a_n}\rangle$, and $g\in G$ where
$$g = x_1^{m_1}\cdots x_n^{m_n}\prod_{1\leq i<j\leq n}
[x_j,x_i]^{m_{ji}}.$$ Then $g\in H$ if and only if  $\forall i,j\in
\{1,\ldots,n\}$, $a_i|m_i$
and $a_ia_j|m_{ji}$.}

\proof Suppose that 
$m_{ji}=a_ja_ik_{ji}$, with $k_{ji}\in\Z$ for each $i$ and $j$. Then
$$g = (x_1^{a_1})^{k_1}\cdots (x_n^{a_n})^{k_n}\prod_{1\leq i<j\leq n}
[x_j^{a_j},x_i^{a_i}]^{k_{ji}},$$
which clearly lies in~$H$.

The converse follows from~(\ref{prodformulainntwo})
and~(\ref{inverseformulaintwo}), after noting that the generators
of~$H$ are of the prescribed~form.\endproof

\thm{domnottrivial}{Let $F(x,y)$ be the relatively free
${\cal N}_2$-group on two generators, and let~$n>1$. Let~$H$ be the
subgroup generated by~$x^n$ and~$y^n$. Then \hbox{$[x,y]^n\in{\rm
dom}_{F(x,y)}^{{\cal N}_2}(H)$}, and $[x,y]^n\notin H$.}

\proof That $[x,y]^n\notin H$ follows from \ref{description}.
Let $G\in{\cal N}_2$, and let
\hbox{$f,g:F(x,y)\rightarrow G$} be two morphisms such that $f|_H =
g|_H$. Then:
$$\eqalign{f\left([x,y]^n\right) &= f\left([x^n,y]\right) = [f(x^n),f(y)]
                                = [g(x^n),f(y)]{\qquad\hbox{(since
$x^n\in H$)}}\cr
                                &= [g(x)^n,f(y)]
                                = [g(x),f(y)]^n\cr}$$
and by a symmetric argument, this term equals $g([x,y]^n)$. In
particular, $[x,y]^n$ lies in ${\rm dom}_{F(x,y)}^{{\cal N}_2}(H)$.\endproof

A similar argument, noting that commutators are central and the image
of a commutator is a commutator, yields:

\cor{domlemma}{Let $G\in{\cal N}_2$, and let $H$ be a subgroup of
$G$. Let $x$ and $y$ be elements of~$G$, $q>0$, and $x',y'\in G'$ such
that $x^qx',y^qy'\in H$. Then
$[x,y]^q$ lies in ${\rm dom}_G^{{\cal N}_2}(H)$.\noproof}

We can get more information thanks to a theorem of Maier, which we 
now~state:

\lemma{niltwoamalg}{{\rm (B.~Maier, Satz 3 in~{\bf [\cite{amalgtwo}]})}
Let $G,K \in{\cal N}_2$, and suppose that $[G,K;H]$ is an amalgam. The
amalgam is strongly embeddable in a group~$M\in {\cal N}_2$ if and
only if the following conditions hold:}
{\it\parindent=30pt
\item{\rm(\numbeq{domeqone})}$G'\cap H \subseteq Z(K)$ and $K'\cap H\subseteq 
Z(G)$.
\item{\rm(\numbeq{domeqtwo})}$\forall q\geq 0$, $g\in G$, $g'\in G'$,
$k\in K$, $k'\in K'$, if
$g^qg', k^qk'\in H$ then both the element $[g,k^qk']$ of $G$ and
the element $[g^qg',k]$ of~$K$ belong to~$H$ and are
equal.\noproof\par}

In fact, condition (\ref{domeqone}) is a consequence of
(\ref{domeqtwo}), as noted by Maier 
(see Remark~3 in~{\bf [\cite{amalgtwo}]}). Set $q=0$, and take
$g=k'=e$ and $g'\in G'\cap H$. Taking any $k\in K$, the hypothesis
of~(\ref{domeqtwo}) is satisfied, 
so we conclude that $[g',k]=[e,e]=e$; so $g'$ commutes with $k$. If we
let~$k$ range over all of~$K$, we have that $g'$ lies in the center of~$K$, 
giving the first inclusion in~(\ref{domeqone}). The second inclusion is proved
in the same way.

\cor{specamal}{Let $G\in {\cal N}_2$ and let~$H$ be a subgroup
of~$G$. Then the special amalgam $[G,G;H]$ is strongly embeddable in 
a group~$M\in{\cal N}_2$ if and only if
$$\forall q\geq 0,\ \forall x,y\in G,\ \forall x',y'\in G',\hbox{\ if\ }
x^qx',y^qy' \in H\hbox{\ then\ } [x,y]^q\in H.\leqno(\numbeq{domclosed})$$}

\proof We will show that if a subgroup~$H$ of~$G$ satisfies
(\ref{domclosed}), then $G$ and~$H$ satisfies
(\ref{domeqtwo}), with~$K=G$. To establish this, we
apply bilinearity of the commutator (we could not do that in the
context of \ref{niltwoamalg}
because although $g^qg'$ lies in $K$, the factors $g^q$ and $g'$ may
not lie there) to get that both $[x^qx',y]$ and $[x,y^qy']$ are equal
to $[x,y]^q$; so condition (\ref{domeqtwo}) becomes just the condition
that the latter element lie in $H$, and this is
condition~(\ref{domclosed}); condition~(\ref{domeqone}) now follows
from~(\ref{domeqtwo}), as noted~above.\endproof

\thm{dominion}{Let $G\in {\cal N}_2$, and let $H$ a subgroup of~$G$. Let $D$
be the subgroup of~$G$ generated by all elements of~$H$ and all
elements of the form $[x,y]^q$, where $x,y\in G$, $q\geq 0$, and there
exists $x',y'\in [G,G]$ such that $x^qx',y^qy'\in H$. Then $D={\rm
dom}_G^{{\cal N}_2}(H).$}

{\it Proof:} It suffices to show that $D$ and $G$
satisfy~(\ref{domclosed}), as then $D$ will be its own dominion, and
it already contains~$H$ and is contained in the dominion of~$H$.

Let $q>0$, and let $x,y\in G$, $x',y'\in [G,G]$ such that $x^qx',y^qy'\in
D$. We want to show that $[x,y]^q\in D$. By \ref{domlemma}, $D\subseteq{\rm
dom}_G^{{\cal N}_2}(H)$, and it follows from
\ref{normalclosure} that we can
perturb $x^qx',y^qy'\in D$ by elements of~$[H,G]$ to get elements
of~$H$; i{.}e{.}, we have $x^qx'', y^qy''\in H$ for appropriate
$x'',y''\in [G,G]$. Hence $[x,y]^q\in D$ by~construction.\endproof

We also note the following property:

\cor{somepowerinsubgroup}{Let $G\in {\cal N}_2$ and let $H$ be a
subgroup of~$G$. If $d$ lies in ${\rm dom}_G^{{\cal N}_2}(H)$, then there
exists $n>0$ such that $d^{\,n}\in H$.}

\proof Let $d\in {\rm dom}_G^{{\cal N}_2}(H)$. By \ref{dominion} there
exists $h\in H$ and 
$$x_1,\ldots,x_m, y_1,\ldots,y_m\in G,\;
x_1',\ldots,x_m',y_1',\ldots,y_m'\in [G,G],\;
q_1,\ldots,q_m>0$$
such that $x_i^{q_i}x_i',y_i^{q_i}y_i'\in H$ for each~$i$, and
$$d = h [x_1,y_1]^{q_1}[x_2,y_2]^{q_2}\cdots[x_m,y_m]^{q_m}.$$
Note that $[x_i,y_i]^{q_i^2}\in H$ for each~$i$. Let $n=\prod
q_i$. Since commutators are central, we have
$$\eqalign{d^{\,n} &= \left(h
[x_1,y_1]^{q_1}[x_2,y_2]^{q_2}\cdots[x_m,y_m]^{q_m}\right)^n\cr
&= h^n[x_1,y_1]^{q_1n}[x_2,y_2]^{q_2n}\cdots[x_m,y_m]^{q_mn}\cr}$$
Since $q_i^2 | q_in$ for each~$i$, it follows that every term on the
right hand side lies in~$H$, and therefore, $d^{\,n}\in H$,
as~claimed.\endproof

Recall that a subgroup of a finitely generated nilpotent group is
finitely generated, and that a finitely generated nilpotent group in
which every element is torsion is~finite.

\cor{finiteindex}{Let $G\in {\cal N}_2$ be any group,
and let~$H$ be any subgroup of~$G$. Then $H$ is normal in~${\rm
dom}_G^{{\cal N}_2}(H)$. If $G$ is finitely generated, then $H$ is of
finite index in~${\rm dom}_G^{{\cal N}_2}(H)$.}

\proof Since ${\rm dom}_G^{{\cal N}_2}(H)$ is generated by~$H$ and
central elements of~$G$, it follows that it normalizes~$H$. This
proves the first assertion. To establish the second assertion, note
that if~$G$ is finitely generated, then so is the dominion of~$H$, and
that ${\rm dom}_G^{{\cal N}_2}(H)/H$ is a torsion group
by~\ref{somepowerinsubgroup}. Therefore, it is a finitely generated
torsion nilpotent group, hence~finite.\endproof
 
{\bf Example \numbeq{ubiquity}.} Our ubiquitous example. Let
$F=F(x,y)$ be the relatively free ${\cal N}_2$-group of rank 2,
and let $n$ be an integer greater than 1. Let~$H$ be the subgroup
of~$F$ generated by $x^n$ and $y^n$. From \ref{description} it follows
that~$H$ consists exactly of all elements of the~form
$$x^{an} y^{bn} [x,y]^{cn^2}\qquad a,b,c\in\Z.$$ Note that this is the
situation we had in \ref{domnottrivial}; there we proved that the
dominion contained the element $[x,y]^n$, which is not in~$H$. Using
\ref{dominion} we can now see that the dominion is actually {\it
generated} by~$H$ and~$[x,y]^n$, which is a subgroup of~$G$ strictly
larger than~$H$. We will use variations of this example
below.\endproof

In \ref{dirproduct} we noted that the dominion construction respects
finite direct products, and added a caution that
the analogous result does not hold for an infinite number of direct
factors. We will now provide the example we promised in
\ref{notinfinite}.

{\bf Example \numbeq{counterexone}.} Let $F$ be the relatively free
${\cal N}_2$ group of rank two generated by $x$ and~$y$.
 
Let $G=F^{\ZZ_+}$ be the direct product of countably many copies
of~$F$, indexed by the positive integers. For each $i>0$, let $H_i$ be
the subgroup of~$F$ generated by $x^i$ and $y^i$. Let $H$ be the
subgroup of~$G$ given by $\prod H_i$.
 
Example~\ref{ubiquity} shows that ${\rm dom}_F^{{\cal N}_2}(H_i) =
\langle x^i,y^i,[x,y]^i\rangle$, that is all elements of~$F$ which can
be written in the form $x^{ai}y^{bi}[x,y]^{ci}$ with $a$, $b$, and~$c$
in~$\Z$.
 
We claim that $${\rm dom}_G^{{\cal N}_2}(H) \propcont \prod_{i>0}{\rm
dom}_F^{{\cal N}_2}(H_i).$$
Consider the element $d=([x,y]^i)_{i>0}$. This is an
element of
$\prod{\rm dom}_F^{{\cal N}_2}(H_i)$ by the discussion in the
preceding paragraph.
To reach a contradiction, assume that $d\in{\rm
dom}_G^{{\cal N}_2}(H)$. From \ref{somepowerinsubgroup} it follows
that there exist some~$n>0$ with $d^{\,n}\in H$. Since 
$d^{\,n} = ([x,y]^{in})_{i>0}$, if $d^{\,n}\in H$, it follows that $i^2|in$
for all~$i>0$, which is clearly~impossible.\endproof

\Section{Dominions in the variety of $2$-Engel groups}{domstwoengel}

The proof that there are instances of nontrivial dominions in the variety
${\cal N}_2$ only relies on the fact that for a group $G\in{\cal N}_2$,
$$\forall x,y\in G\;\forall n\in\hbox{\Z}\qquad
[x^n,y]=[x,y]^n=[x,y^n];\leqno(\numbeq{engellawtwo})$$
so we might ask which are the groups that satisfy~(\ref{engellawtwo}).

\lemma{ifforminusforall}{Let $G$ be a group that satisfies
$$\forall x,y\in G\qquad \left[x^{-1},y\right]=\left[x,y\right]^{-1} =
\left [x,y^{-1}\right].$$
Then $G$ satisfies the identity $[[x,y],y]=e$. In particular, if~$G$
satisfies~{\rm (\ref{engellawtwo})}, it also satisfies the identity
$[[x,y],y]=e$.}

\proof In any group we have that $[x,y]=[y,x]^{-1}$, and by
\ref{commident}(b) and our hypothesis we~have
$[x,y]=[y,x]^{-1} = [y^{-1},x]$. Therefore,
$$\eqalignno{[x,y]^y &= y^{-1}[x,y]y
                     = y^{-1}[y^{-1},x]y
                     = y^{-1}yx^{-1}y^{-1}xy\cr
                     &= x^{-1}y^{-1}xy
                     = [x,y].\cr}$$
Hence $[x,y]$ commutes with~$y$, so $[[x,y],y]=e$ as~claimed.\endproof

\defn A group $G$ is said to be {\it 2-Engel\/} if it satisfies the law
$$\left[\left[x,y\right],y\right]=e.$$
In general, the $k$-th Engel Law is given by
$[x,y,\ldots,y]=e$, where the $y$ occurs $k$~times.

In fact, we have a converse to \ref{ifforminusforall}:

\lemma{engel}{{\rm (Lemma~5.42 in~{\bf [\cite{rotman}]})} Let $G$ be a group in
which $[x,y]$
commutes with~$y$. Then $\forall n\in\hbox{\Z}$, $[x,y^n]=[x,y]^n$.\noproof}

\cor{engelident}{A group $G$ is 2-Engel if and only if for all $x$,
$y\in G$, and for all $n\in\hbox{\Z}$, \hbox{$[x,y^n]=[x,y]^n=[x^n,y]$}.\noproof}

The reader may find it amusing to verify the
following~result:

\lemma{tangential}{For a group~$G$ the following are equivalent:}
\nobreak
{\parindent=30pt\it
\item{(i)}{$G$ is a 2-Engel group.}\par
\item{(ii)}{$G$ satisfies~{\rm (\ref{engellawtwo})}.}\par
\item{(iii)}{For all $x,y\in G$,
$[x,y^{-1}]=[x,y]^{-1}=[x^{-1},y]$.}\par
\item{(iv)}{For all $x,y\in G$, $[x,y^2]=[x,y]^2=[x^2,y]$.}\par
\item{(v)}{There exists an integer $n$ such that for all $x,y\in G$,
and \hbox{$i=0,1,2$}, $$[x,y^{n+i}]=[x,y]^{n+i} =
[x^{n+i},y].$$ (Cf.~Exercise~2.3.4, pp.~31
in~{\bf [\cite{herstein}]})}\par
\item{(vi)}{For all $x,y\in G$, $[x^{-1},y]=[x,y^{-1}]$.}\par
\item{(vii)}{For all $x,y\in G$, $[x^2,y]=[x,y^2]$.}\par
\item{(viii)}{Every subgroup of~$G$ which is generated by at most two
elements is nilpotent of class at most~$2$.\noproof}\par}

From the comments above, it follows that:

\thm{domsnontrivinengel}{The variety of $2$-Engel groups has instances
of nontrivial dominions.\noproof}

\Section{Dominions of subgroups of ${\cal N}_2$-groups in ${\cal A}^2\cap{\cal
N}_c$}{domsniltwoinmetab}

We turn our attention to the variety of metabelian nilpotent groups of
class at most~$c$, with $c>1$. In this section we will prove that
these varieties also contain instances of nontrivial dominions. We
will in fact prove a bit more: that there is a finitely generated
group~$G$, nilpotent of class~$2$, such that for each given $c>1$,
there exists a subgroup~$H_c$ of~$G$ such that
$$H_c\propcont {\rm dom}_G^{{\cal N}_2}(H_c) = {\rm dom}_G^{{\cal A}^2\cap
{\cal N}_c}(H_c).$$
 
In fact, we will prove the corresponding result with the variety ${\cal
A}^2\cap {\cal N}_c$ replaced by the variety~${\cal N}_c$
itself in \ref{domsniltwoinnilk}. This will of course imply it for the case
we are now contemplating. However, the technique and calculations are
more transparent in the metabelian case. Thus we include it
first as an introduction to the method, rather than deducing it from
the more general result~later.

\label{Definition\noexpand~\the\sectno.\the\dispno}{defforref}
\defn Let $G$ be a group, and $x$ and $y$ elements of~$G$. We~write
$$[x,{}_ny]= [x,\underbrace{y,\ldots,y}_{n\ times}].$$

In a metabelian group, of all the basic commutators only the left
normed may be nontrivial. The left normed basic commutators on $x_1$,
$x_2,\ldots,x_n$ will look as follows:
$$[x_{i_1},x_{i_2},\ldots,x_{i_m}]$$
where $i_1>i_2$, and $i_2\leq i_3\leq i_4\leq\cdots\leq i_m$.

The following lemma is easy to establish using \ref{commident}:

\lemma{commidentmetabelian}{Let $G$ be a group satisfying
$[G_2,G_3]=\{e\}$, and let~$x$, $y$, $z$, $w$, and~$t$ be elements
of~$G$. Then
$$\eqalign{[[x,y][z,w],t] &= [x,y,t][z,w,t],\cr [t,[x,y][z,w]] &=
[t,[x,y]][t,[z,w]].\cr}$$ In particular, the above identities hold
in any metabelian~group.\noproof}
 
\lemma{inductivecomm}{Let $G$ be a group satisfying $[G_2,G_3]=\{e\}$, and
$x_1,\ldots,x_m$ be elements of~$G$, with $m\geq 3$. Then for any
$n>0$ the following identities hold:
$$\eqalignno{[[\cdots[[x_1,x_2],x_3]\ldots],x_m]^n &=
[[\cdots[[x_1,x_2],x_3]\ldots]^n,x_m]\cr
&\vdots\cr
&=[[\cdots[[x_1,x_2]^n,x_3]\ldots],x_m].&\hbox{\noproof}\cr}$$}

The following result is easily proven using induction:
 
\lemma{identmetabelian}{Let $G$ be a group such that $[G_2,G_3]=\{e\}$,
and let $x$
and~$y$ be elements of~$G$. Then for any $n>0$, the following
identities~hold:
$$\leqalignno{[y^n,x] &= [y,x]^{n\choose 1}[y,x,y]^{n\choose
2}[y,x,{}_2y]^{n \choose 3}\cdots[y,x,{}_{(n-1)}y]^{n\choose
n};&(\numbeq{metaxncommyprime})\cr [y,x^n] &= [y,x]^{n\choose
1}[y,x,x]^{n\choose 2}[y,{}_3x]^{n\choose 3}\cdots [y,{}_nx]^{n\choose
n}.&(\numbeq{metaxcommynprime})\cr}$$ In particular,
{\rm (\ref{metaxncommyprime})} and {\rm(\ref{metaxcommynprime})} hold in any
metabelian~group.\noproof}
 
Note that (\ref{metaxncommyprime}) and~(\ref{metaxcommynprime}), in a
group satisfying $[G_2,G_3]=\{e\}$ may be rewritten
$$\leqalignno{[y,x]^n &= [y^n,x][y,x,y]^{-{n\choose
2}}[y,x,{}_2y]^{-{n \choose 3}}\cdots[y,x,{}_{(n-1)}y]^{-{n\choose
n}};&(\numbeq{metaxncommy})\cr 
[y,x]^n &= [y,x^n][y,x,x]^{-{n\choose 2}}[y,{}_3x]^{-{n\choose
3}}\cdots [y,{}_nx]^{-{n\choose n}}.&(\numbeq{metaxcommyn})\cr}$$

The variety ${\cal A}^2$ mentioned in all further results in this
section could be replaced with the larger variety of all groups~$G$ in
which $[G_2,G_3]=\{e\}$. We state the results for~${\cal A}^2$
for~simplicity. 

\lemma{thirdtermexppmetab}{Let $c>1$ and let $G\in{\cal A}^2\cap {\cal
N}_c$. Let $y,z\in G$, $p$ a prime such that $p\geq c$. If
$$[z,y^p,z]= [z^p,y,y]=e\leqno(\numbeq{condonzandywithp})$$
then $\langle y,z\rangle_3$ has exponent~$p$.}

\proof Since $G$ is metabelian, and $\langle y,z\rangle_3$ is abelian
and generated by the basic commutators in~$y$ and~$z$ of ${\rm weight}\geq
3$, it suffices to verify that all basic commutators of ${\rm weight}\geq 3$
are of exponent~$p$. This is vacuous for all commutators of weight
greater than~$c$, regardless of whether~(\ref{condonzandywithp})
holds. Suppose the lemma is false, and let $m_0\in\Z$,
with $3\leq m_0\leq c$ be maximal with respect to the property that
there exists a basic commutator~$c_0$ in~$y$ and~$z$ of weight~$m_0$ which
is not of exponent~$p$.

As noted following \ref{defforref}, $c_0=[z,{}_ny,{}_kz]$; since ${\rm
weight}(c_0)=m_0$, we must~have $n+k=m_0-1\geq 2$. Therefore, if
$n>1$ we have
$$\eqalignno{c_0^{\,p} &= [z,{}_ny,{}_kz]^p
               = \bigl[[z,y]^p,{}_{(n-1)}y,{}_kz\bigr]\quad\hbox{(by
                   \ref{inductivecomm})}\cr
               &= \bigl[[z^p,y][z,y,z]^{-{p\choose
                   2}}\cdots[z,y,{}_{(c-2)}z]^{-{p\choose
                   c-1}},{}_{(n-1)}y,{}_kz\bigr]\quad\bigl(\hbox{by
(\ref{metaxcommyn})}\bigr)\cr 
               &= \bigr[[z^p,y],{}_{(n-1)}y,{}_{k}z\bigr]\cdots
               \bigr[[z,y,{}_{(c-2)}z]^{-{p\choose
                   c-1}},{}_{(n-1)}y,{}_kz\bigr]\cr
               &=[z^p,{}_ny,{}_kz][z,y,z,{}_{(n-1)}y,{}_kz]^{-{p\choose
               2}}\cdots[z,y,{}_{(c-2)}z,{}_{(n-1)}y,{}_kz]^{-{p\choose 
               c-1}}\cr 
&\qquad\bigl(\hbox{by \ref{inductivecomm}}\bigr).\cr}$$

The first factor is trivial, since $[z^p,y,y]=e$. All the subsequent
factors are commutators in~$y$ and~$z$ of ${\rm weight}>m_0$. In
particular, they are all of exponent~$p$ by choice of~$m_0$. Since
$p\geq c$, the powers $-{p\choose i}$ to which they are raised are
multiples of~$p$, hence all the factors on the right hand side are
trivial. Thus, $c_0^{\,p}=e$. If, on the other hand, $n=1$, then we
simply use~(\ref{metaxncommy}) instead of~(\ref{metaxcommyn}) to
``pull in'' the exponent, and proceed as~above, to conclude that
$c_0^{\,p}=e$, a contradiction to the choice of $c_0$.\endproof

\rmrk{othertriples} It is worth noting that \ref{thirdtermexppmetab}
also holds if~(\ref{condonzandywithp}) is replaced by any one of:
$$[z,y^p,z]=[z,y^p,y]=e;\; [z^p,y,z]=[z^p,y,y]=e;\;
[z^p,y,z]=[z,y^p,z]=e.$$

The only difference in the proof would be that
under the first condition, one would handle both the $n>1$ and
the $n=1$ case as the $n>1$ case is handled in
\ref{thirdtermexppmetab}; under the second one, we would
handle them both like the $n=1$ case; and
under~the last condition, one would interchange the methods used
for $n>1$ and $n=1$. However, only~(\ref{condonzandywithp}) is
relevant to our study of dominions in~${\cal A}^2\cap{\cal N}_c$.
In any case, all four conditions are equivalent, which will follow
from \ref{thirdtermexppmetab} and \ref{struiklemhtwo}.

\cor{fromptoengelcond}{Let $c>1$ and let $G\in {\cal A}^2\cap {\cal
N}_c$. Let $y,z\in G$, $p$ a prime such that $p\geq
c$. If~{\rm (\ref{condonzandywithp})} holds, then
$[z,y^p]=[z,y]^p=[z^p,y]$.}

\proof From~(\ref{metaxncommyprime})
and~(\ref{metaxcommynprime}) we have
$$\eqalign{[z^p,y]&= [z,y]^p [z,y,z]^{p\choose 2}[z,y,z,z]^{p\choose
3}\cdots
[z,y,{}_{(c-2)}z]^{p\choose c-1}\cr
[z,y^p]&= [z,y]^p [z,y,y]^{p\choose 2}[z,y,y,y]^{p\choose 3} \cdots
[z,y,{}_{(c-2)}y]^{p\choose c-1}.\cr}$$

Every factor on the right hand side of both expressions is
raised to a multiple of~$p$, and that all but the first term lie
in~$\langle z,y\rangle_3$. Since the latter subgroup is of
exponent~$p$ by \ref{thirdtermexppmetab}, it follows that
$[z,y^p]=[z,y]^p=[z^p,y]$
as~claimed.\endproof

\thm{dommetanontrivial}{Fix $c>1$. Let $p$ be a prime, $p\geq
c$, and let $G=F^{\scriptscriptstyle (2)}(x,y)$ be the relatively free
nilpotent group of class~$2$ on
two generators. Let $H$ be the subgroup of~$G$ generated by~$x^p$
and~$y^p$. Then ${\rm dom}_G^{{\cal A}^2\cap {\cal N}_c}(H)$
is generated
by~$x^p$, $y^p$, and~$[x,y]^p$. In particular, it equals ${\rm
dom}_G^{{\cal N}_2}(H)$ and properly contains~$H$.}
 
\proof Let ${\cal V}={\cal A}^2\cap {\cal N}_c$. We look at the
amalgamated coproduct of two copies of $G$ 
over the subgroup~$H$ in the variety ${\cal V}$, denoted by
$F=G\amalg_H^{{\cal V}} G$. This is the group
with presentation
$$\Bigl\langle x,y,z,w\,\Bigm|\, x^p=z^p;\  y^p=w^p;\quad
\bigl(\langle
x,y\rangle\bigr)_2=\{e\};\ \bigl(\langle
z,w\rangle\bigr)_2=\{e\}\Bigr\rangle$$
in the variety determined by the identities
$$[[x_1,x_2],[x_3,x_4]]=e\quad {\rm and} \quad
[x_1,\ldots,x_c,x_{c+1}]=e;$$ 
where by
$\bigl(\langle x,y\rangle\bigr)_2=\{e\}$ we mean that the
subgroup generated by $x$ and~$y$ is nilpotent of class two,~etc.
 
Then $[x,y]^p\in {\rm dom}_{G}^{\cal V}({H})$ if and only if
$[x,y]^p=[z,w]^p$ in
$F$. We note that, since the subgroup generated by~$x$ and~$y$, and
that generated by~$z$ and~$w$, are both nilpotent of class two, if
follows~that
$$[x,y]^p = [x^p,y] = [z^p,y],\qquad
[z,w]^p = [z,w^p]= [z,y^p]$$
so we want to see whether $[z^p,y]=[z,y^p]$.
Note that $$[z,y^p,z]=[z,w^p,z]=e\quad\hbox{and}\quad
[z^p,y,y]=[x^p,y,y]=e.$$  
In particular, by \ref{thirdtermexppmetab}, $\langle z,y\rangle_3$ is of
exponent~$p$, and by \ref{fromptoengelcond}, $[z^p,y]=[z,y^p]$. 
Hence $[x,y]^p\in {\rm dom}_G^{{\cal A}^2\cap {\cal
N}_c}(H)$. 

Therefore, by Example~\ref{ubiquity}, 
${\rm dom}_G^{{\cal N}_2}(H)\subseteq {\rm dom}_G^{{\cal A}^2\cap
{\cal N}_c}(H)$; the reverse inclusion now follows from
\ref{varinclusions}.\endproof

Note that, if we take $G=F^{(2)}(x,y)$ and $H=\langle x^p,y^p\rangle$,
and we let $c>p$, then the argument in \ref{thirdtermexppmetab} breaks
down. This raises the question of whether for any $G\in {\cal N}_2$,
any subgroup~$H$ of~$G$, and any given $g$ in $G\backslash
H$, there exists a $k$ (depending on $g$) such that $g$ is not in the
dominion of $H$ in~$G$, in the variety ${\cal A}^2\cap {\cal
N}_{k}$. Put another way, we ask whether dominions of subgroups of
${\cal N}_2$ groups are trivial in the category of all
met\discretionary{a-}{a}{a}belian nilpotent~groups (note that this
category is not a variety). We will answer this question in \ref{nontrivnil}.

\Section{Dominions in ${\cal A}^2\cap {\cal N}\!il$}{domsmetabnil}

In this section we will investigate dominions in the category of all
metabelian nilpotent groups. We take the oportunity to obtain some
results about dominions in the variety of all metabelian groups as well.
We begin by making some general observations on ${\cal A}^2$:

\lemma{metabcomm}{{\rm(Lemma~34.51 in~{\bf [\cite{hneumann}]})} Let
$G$ be metabelian, and let $x\in G'$. Then for all
elements $y_1,\ldots,y_n\in G$ and every permutation of~$n$
elements~$\sigma\in S_n$,
$$[x,y_1,\ldots,y_n] =
[x,y_{\sigma(1)},\ldots,y_{\sigma(n)}].\eqno{\noproof}$$} 

\lemma{generate}{{\rm (Theorem~36.33 in~{\bf [\cite{hneumann}]})} For
a relatively free group $F$ in the variety ${\cal A}^2$,
the left normed basic commutators of ${\rm weight}\geq 2$ freely
generate a free abelian subgroup of the derived group~$F'$.\noproof}
 
We add a caution:
 
\lemma{notallbasiccomms}{{\rm (Proposition~36.24 in~{\bf
[\cite{hneumann}]})} The
basic commutators of weight at least two in~$F$, the free group of rank
two in ${\cal A}^2$, do not generate the derived group~$F'$, even
though they generate $F'$ modulo every term of the lower
central~series.\noproof}

We will first obtain results about dominions in ${\cal A}^2$.
We will then derive the results we want by
restricting our attention to groups within ${\cal A}^2\cap {\cal
N}\!il$.

\lemma{domsinmetab}{Let $G \in {\cal A}^2$, $y\in G'$, $x,z\in G$, and
let
$H$ be a subgroup of~$G$. If $y$, $[y,x]$ and~$[y,z]$ lie in~$H$, then
$[y,x,z]\in {\rm dom}_G^{{\cal A}^2}(H)$.}
 
\proof Let $K\in {\cal A}^2$, and let $f,g\colon G\to K$ be two group
morphisms such that $f|_H = g|_H$. Then
$$\eqalignno{f\Bigl([y,x,z]\Bigr)&=
\Bigl[f\bigl([y,x]\bigr),f(z)\Bigr]
=\Bigl[g\bigl([y,x]\bigr),f(z)\Bigr]\qquad \hbox{(since $[y,x]\in
H$)}\cr
&=\Bigl[g(y),g(x),f(z)\Bigr]
=\Bigl[f(y),f(z),g(x)\Bigr]\cr
&\qquad\qquad \hbox{(using \ref{metabcomm} and
the fact that $y\in H$)}\cr
&=\Bigl[f\bigl([y,z]\bigr),g(x)\Bigr]
=\Bigl[g\bigl([y,z]\bigr),g(x)\Bigr]\quad \hbox{(since
$[y,z]\in
H$)}\cr
&=g\Bigl([y,z,x]\Bigr)
=g\Bigl([y,x,z]\Bigr)\quad \hbox{(by
\ref{metabcomm})}\cr}$$
so $[y,x,z]\in {\rm dom}_G^{{\cal A}^2}(H)$, as claimed.\endproof
 
\thm{domsinmetabgen}{Let $G\in {\cal A}^2$, $x\in G'$, and
$y_1,\ldots,y_n\in G$. Let
$$H=\Bigl\langle x, [x,y_1], \ldots, [x,y_n]\Bigr\rangle.$$
Then $[x,y_1,\ldots,y_n]\in {\rm dom}_G^{{\cal A}^2}(H)$.}
 
\proof Follows from \ref{domsinmetab} by induction on~$n$.\endproof
 
The next two results express essentially the same thing as
\ref{domsinmetab} and \ref{domsinmetabgen}, but in terms of
conjugation rather than commutation.
 
\lemma{domsinmetabtwo}{Let $G\in {\cal A}^2$, and let $H$ be a
subgroup
of~$G$. Let $x\in G'$, and let $y$ and~$z$ be elements of~$G$. If $x$, $x^y$,
and~$x^z$ lie in~$H$, then $x^{yz}$ lies in ${\rm dom}_G^{{\cal
A}^2}(H)$.}
 
\proof Recall that $x^y=y^{-1}xy=x[x,y]$. Applying
\ref{commident}(d), we obtain
$$x^{yz}=x[x,z][x,y][x,y,z].$$
 
Since $x$ and $x^y$ both lie in~$H$, it follows that $[x,y]$ also lies
in~$H$. Similarily, $[x,z]$ lies in~$H$. By \ref{domsinmetab},
$[x,y,z]\in
{\rm dom}_G^{{\cal A}^2}(H)$, so $x^{yz}\in {\rm dom}_G^{{\cal
A}^2}(H)$, as~claimed.\endproof
 
\thm{domsinmetabtwogen}{Let $G\in {\cal A}^2$ and let $H$ be a
subgroup
of~$G$. Let $x\in G'$, and $y_1,\ldots,y_n$ be elements of~$G$. If
$x\in H$ and $x^{y_i}\in H$ for $i=1,\ldots,n$, then $x^{y_1\cdots
y_n}$ also lies in~${\rm dom}_G^{{\cal A}^2}(H)$.\noproof}
 
\thm{nontrivmetab}{Let $G=F(x,y)$ be the free metabelian group on two
generators. Let $H$ be the subgroup of~$G$ generated by $[y,x]$,
$[y,x,y]$ and $[y,x,x]$. Then $H\propcont {\rm dom}_G^{{\cal
A}^2}(H)$. Moreover, ${\rm dom}_G^{{\cal A}^2}(H) = G'$,
so ${\rm dom}_G^{{\cal A}^2}(H)/H$ is an abelian group of
infinite~rank.}
 
\proof By \ref{generate}, $H$ is a free abelian group on the
generators $[y,x]$,
$[y,x,y]$ and~$[y,x,x]$. Applying \ref{domsinmetabgen} to $H$ and
these generators, we see
that $D={\rm dom}_G^{{\cal A}^2}(H)$ contains all basic commutators of
weight $k\geq 2$.
 
We claim that $G'$ is generated by all elements of the form $[y,x]^g$,
where $g$ ranges over all elements of~$G$. 
Indeed, using \ref{commident}(c) and noting that
$[y,x]^{-1}=[x,y]$, we can decompose any
commutator $[z,w]$ into a product of conjugates of $[y,x]$,
$[y,x^{-1}]$, $[y^{-1},x]$, $[y^{-1},x^{-1}]$ and their inverses. But
$[y,x^{-1}] = ([y,x]^{x^{-1}})^{-1}$, and similarily for $[y^{-1},x]$
and $[y^{-1},x^{-1}]$, so that $G'$ is indeed generated by all
elements of the form~$[y,x]^g$.

Since $D$ contains all basic commutators, it certainly contains
$[y,x]^x$ and
$[y,x]^y$. It will suffice to show that it
also contains $[y,x]^{y^{-1}}$ and $[y,x]^{x^{-1}}$.
 
Let $u=[y,x]^y$. Then $u\in D\cap G'$, and $u^{y^{-1}}= [y,x] \in
D$. By \ref{domsinmetabtwo}, the dominion of~$D$ 
must also contain $u^{y^{-1}y^{-1}}=u^{y^{-2}} =
[y,x]^{y^{-1}}$. Since $D$ is its own dominion, this elements lies
in~$D$. A
similar calculation yields $[y,x]^{x^{-1}}\in D$. This proves that
$G'\subseteq D$. Since $G'$ is normal in~$G$, ${\rm dom}_G^{{\cal
A}^2}(G')=G'$, so $D=G'$ as~claimed.
 
By \ref{generate}, it follows that $D$ contains a free abelian group
on
countably many generators, and $H$ is the subgroup generated by three
of
these generators. Therefore, $D/H$ is abelian of infinite rank,
as~claimed.\endproof
 
\cor{nilinmetab}{There are nontrivial dominions in ${\cal
A}^2\cap{\cal N}\!il$. Specifically, if $F$ is the relatively free
${\cal N}_4$ group of rank 2, with generators $x$ and~$y$, then the
dominion of the subgroup generated by $[x,y]$, $[x,y,x]$, and
$[x,y,y]$ is the entire derived subgroup of~$F$.}

\proof $F$ is the quotient of the relatively free ${\cal A}^2$ group
on two generators, modulo the fifth term of the lower central
series. By \ref{quotient}, we get the result.\endproof

\rmrk{observation} Note that this example also establishes that
the finite index clause of \ref{finiteindex} does not hold for general
varieties of nilpotent groups, by looking at ${\cal A}^2\cap{\cal
N}_4$.

\Section{Dominions of subgroups of ${\cal N}_2$-groups in ${\cal
N}_c$}{domsniltwoinnilk}

As promised at the beginning of \ref{domsniltwoinmetab}, we will now
prove the analog of \ref{dommetanontrivial} for the varieties
${\cal N}_c$. The proof is indeed very similar to that of
\ref{dommetanontrivial}, with the added complications that arise from
the nonabelian nature of the commutator subgroup in the more
general~case.

We first require two technical lemmas, which we now state:

\lemma{struikhthree}{{\rm (Struik, Theorem~H3 in~{\bf
[\cite{struikone}]})}
Let $x_1,\ldots,x_s$ be any elements of a nilpotent group~$G$, and let
$u_1<u_2<\ldots$ be a system of basic commutators on the~$x_i$. Let
$\sigma$ be a fixed permutation of $\{1,2,\ldots,s\}$, and let
$n>0$. Then
$$(x_1x_2\cdots x_s)^n = x_{\sigma(1)}^n x_{\sigma(2)}^n\cdots
x_{\sigma(s)}^nu_1^{f_1(n)}u_2^{f_2(n)}\ldots\leqno(\numbeq{powerin})$$
where 
$$f_i(n) = a_1{n\choose 1} + a_2{n\choose 2} +\cdots +
a_{w_i}{n\choose w_i}\leqno(\numbeq{exponentsofcomms})$$
with $a_k$ an integer, and $w_i$ the weight of~$u_i$ in the
$x_j$.\noproof}

\lemma{struiklemhtwo}{{\rm (Struik, Lemma~H2 in~{\bf
[\cite{struikone}]})}
Let $n$ be a fixed integer and let~$G$ be a nilpotent group of
class~$c$. If $x_j\in G$, then
$$[x_1,\ldots,x_{i-1},x_i^n,x_{i+1},\ldots,x_r] = [x_1,\ldots,x_r]^n
v_1^{f_1(n)} v_2^{f_2(n)}\ldots\leqno(\numbeq{poweringen})$$
where the
$v_k$ are basic commutators in $x_1,\ldots,x_r$ of ${\it weight}>r$,
and every $x_j$ appears in each commutator $v_k$ for $1\leq j\leq
r$. The $f_i$ are of the form~{\rm (\ref{exponentsofcomms})} where $w_i$ is
the weight of~$v_i$ minus $(r-1)$.\noproof}

We will find~(\ref{poweringen}) useful in situations when we have
commutators in some terms, some of which are shown as powers, and
we want to ``pull the exponent out.'' But at other times we will want
to reverse this process and pull exponents ``into'' a commutator. In
such situations, we will use~(\ref{poweringen}) to express
$[b_1,\ldots,b_r]^n$ in terms of other commutators. Let us call the
resulting identity~$(\ref{poweringen}')$; that is,
$$[b_1,\ldots,b_r]^n =
[b_1,\ldots,b_{i-1},b_i^n,b_{i+1},\ldots,b_r]\cdots v_2^{-f_2(n)}
v_1^{-f_1(n)}.\leqno(\ref{poweringen}')$$

There is a slight refinement to the above formulas for
the case $r=2$ and $n$ a prime. Namely, 
 
\lemma{struikeq}{{\rm (Struik; see equations labeled (57) and (58)
in~{\bf [\cite{struiktwo}]})} Let $F$ be the absolutely
free group with
generators $x$ and~$y$, and let $p$ be a prime. Then
$$\leqalignno{[y^p,x]&\equiv [y,x]^p\left(\prod_{i=4}^{q(p+1)}
c_i^{p\beta_i}\right) [y,x,{}_{p-1}y]
\pmod{F_{p+2}}&(\numbeq{struikxpcommy})\cr
[y,x^p] &\equiv [y,x]^p\left(\prod_{i=4}^{q(p+1)}
{c_i}^{p\alpha_i}\right) [y,{}_px]
\pmod{F_{p+2}}&(\numbeq{struikxcommyp})\cr}$$
where $q(n)$ is the number of basic
commutators of weight less than or equal to~$n$ in two symbols $x$ and~$y$,
$c_4<c_5<\ldots<c_{q(p+1)}$ are the basic commutators of ${\it
weight}>2$ and ${\it weight}\leq p+1$, and $\alpha_i$ and $\beta_i$
are~integers.\noproof}

To finish our preparatory lemmas, we need a description of the basic
commutators on two generators. This description is easily established
by induction on the weight:

\lemma{freeontwo}{Let $F(x,y)$ be the free group on two
generators. Then every basic commutator of weight$\,\geq 3$ is of the
form
$$[y,x,y,c_4,\ldots,c_r]\quad {\rm or}\quad [y,x,x,c_4,\ldots,c_r]
\leqno(\numbeq{shapeofcomms})$$
where $c_4,\ldots,c_r$ are basic commutators in~$x$ and~$y$.\noproof}

We now proceed as in \ref{domsniltwoinmetab}:

\lemma{thirdtermexpp}{Let $G\in {\cal N}_c$ and let $p$ be a prime,
with $p\geq c$. If $z,y\in G$ are such that
$[z,y^p,z] = [z^p,y,y]=e$,
then $\langle z,y\rangle_3$ is of exponent~$p$.}

\proof Assume the lemma is false. Since every commutator of ${\rm
weight}>c$ has exponent~$p$, any counterexample would have smaller
weight. Let $m_0$ be maximal with respect to the property that there
exist an element $c_0$ in $\langle z,y\rangle_{m_0}$ which is not of
exponent~$p$.

First we prove that a basic commutator $c_0$ of weight exactly~$m_0$ is of
exponent~$p$, by writing it as in~(\ref{shapeofcomms}) and then
using~$(\ref{poweringen}')$ to ``pull in'' the exponent, and express
$c_0^{\,p}$ as a product of a commutator which is trivial by
hypothesis, and commutators of higher weight each raised to a power
which is a multiple of~$p$, to reach a contradiction.

Next we assume that $c_0$ is a product of basic commutators of weight
at least~$m_0$, and apply~(\ref{powerin}) to reach a~contradiction.\endproof

\cor{fromptoengelcondgeneral}{Let $c>0$, and let $G\in {\cal N}_c$,
$p$ a prime with $p\geq c$. If $z,y\in G$ are such that
$$[z,y^p,z]=[z^p,y,y]=e,$$
then $[z,y^p]=[z,y]^p=[z^p,y]$.}

\proof By \ref{thirdtermexpp}, $\langle z,y\rangle_3$ is of exponent~$p$.
We now apply Struik's formulas from \ref{struikeq}, to
get
$$[z^p,y] = [z,y]^p
\left(\prod_{i=4}^{q(c)}c_i^{p\beta_i}\right),\qquad
[z,y^p] = [z,y]^p\left(\prod_{i=4}^{q(c)}c_i^{p\alpha_i}\right).$$

Each $c_i$ is of weight at least three, hence of
exponent~$p$. Therefore, we conclude that
$[z,y^p]=[z,y]^p=[z^p,y]$ as~claimed.\endproof

\thm{domnsubc}{Fix $c>1$. Let $p$ be a prime, $p\geq
c$, and let $G=F^{\scriptscriptstyle (2)}(x,y)$ be the relatively free
nilpotent group of class~two on
two generators. Let $H$ be the subgroup of $G$ generated by $x^p$ and
$y^p$. Then ${\rm dom}_{G}^{{\cal N}_c}(H)$ is generated
by $x^p$, $y^p$, and $[x,y]^p$. In particular, it equals
${\rm dom}_G^{{\cal N}_2}(H)$, and properly contains~$H$.}
 
\proof Once again, we look at the amalgamated coproduct of two copies
of $G$
over the subgroup~$H$ in the variety ${\cal N}_c$, denoted by
$F=G\amalg_H^{{\cal N}_c} G$. This is the group
with presentation
$$\Bigl\langle x,y,z,w\,\Bigm|\, x^p=z^p;\  y^p=w^p;\quad
\bigl(\langle
x,y\rangle\bigr)_2=\{e\};\ \bigl(\langle
z,w\rangle\bigr)_2=\{e\}\Bigr\rangle$$
in the variety determined by the identity
$[x_1,\ldots,x_c,x_{c+1}]=e$;
where by $\bigl(\langle x,y\rangle\bigr)_2=\{e\}$ we again mean that the
subgroup generated by $x$ and~$y$ is nilpotent of class~2,~etc.

Then $[x,y]^p\in {\rm dom}_{G}^{{\cal N}_c}({H})$ if and only if
$[x,y]^p=[z,w]^p$ in $F$. The rest of the proof now proceeds like the
proof of \ref{dommetanontrivial}.\endproof

\cor{lotsofvarsbetween}{For every $c>1$, any category ${\cal V}$ with
${\cal N}_2\subseteq {\cal V}\subseteq {\cal N}_c$ has instances of
nontrivial dominion. Namely, if $G=F^{\scriptstyle (2)}(x,y)$ is the
relatively free ${\cal N}_2$ group on two generators, and $p$ is a
prime with $p\geq c$, then
$$\langle x^p,y^p\rangle \propcont {\rm dom}_G^{{\cal
N}_2}\bigr(\langle x^p,y^p\rangle\bigr) = {\rm dom}_G^{\cal
V}\bigr(\langle x^p,y^p\rangle\bigr).\eqno\noproof$$}

\rmrk{solvablermrk} \ref{lotsofvarsbetween} includes, among others,
the varieties ${\cal A}^k\cap{\cal N}_c$ of all groups
which are nilpotent of class at most~$c$ and solvable of length at
most~$k$.

The same comments as in the closing of \ref{domsniltwoinmetab} apply
here. So we ask whether dominions of subgroups of ${\cal N}_2$ groups
are trivial in the category of all nilpotent groups.
We will answer this question in
\ref{nontrivnil}.

\Section{Dominions of subgroups of ${\cal N}_2$-groups in ${\cal
N}\!il$}{nontrivnil}

In this section we will prove results similar to
\ref{domnsubc} without the assumption that $p\geq
c$. The arguments are very similar to those used above, and so we will
again only sketch the proofs.

\lemma{generaltermexp}{Let $c>1$, $p$ a prime, and let $a>0$. Let $G\in
{\cal N}_c$, and $y,z\in G$. Suppose that for every $i\geq a$,
$$[z^{p^i},y,y] = [z,y^{p^i},z]=e.\leqno(\numbeq{commsthreetriv})$$
Let $m_0\geq 3$. If $\langle z,y\rangle_{m_0+1}$ is
of exponent $p^k$, and $k\geq a$, then $\langle z,y\rangle_{m_0}$ is
of exponent $p^{{\rm ord}_p(c!)+k}$.}

\proof The first step is to show that every basic commutator of weight
exactly $m_0$ has exponent $p^{{\rm ord}_p(c!)+k}$. We do this by
applying $(\ref{poweringen}')$ in the same manner as we did before,
using the term ${\rm ord}_p(c!)$ to make sure we can factor out at
least a $p^k$ from the exponent $f_i(n)$.
Then, to prove that an arbitrary element of $\langle z,y\rangle_{m_0}$
is also of exponent $p^{{\rm ord}_p(c!)+k}$, we write it as a product
of basic commutators of weight at least~$m_0$ and
apply~(\ref{powerin}) and the fact that the result holds for basic
commutators of weight at least~$m_0$.\endproof

\thm{expofthirdterm}{Let $c>1$, $G\in {\cal N}_c$, $z,y\in G$,
$a>0$ and $p$ a prime. Suppose that for each $i\geq a$, $G$
satisfies~{\rm (\ref{commsthreetriv})}. 
Then $\langle z,y\rangle_3$ is of exponent $p^{(c-3){\rm
ord}_p(c!)+a}$.}

\proof Since $G\in {\cal N}_c$, we have $\langle z,y\rangle_c = \{e\}$,
so it is of exponent~$1$. However, in order to apply
\ref{generaltermexp} we note that it is also of exponent $p^a$.

By \ref{generaltermexp}, $\langle z,y\rangle_{c-1}$ is of exponent
$p^{{\rm ord}_p(c!)+a}$; $\langle z,y\rangle_{c-2}$ is of exponent
$p^{2{\rm ord}_p(c!)+2}$; etc. Continuing this until the third term,
we obtain that $\langle z,y\rangle_3$ is of exponent
$$p^{(c-3){\rm ord}_p(c!)+a}$$
as~claimed.\endproof

\cor{transferifenough}{Let $c>1$, $G\in {\cal N}_c$, $z,y\in G$,
$a>0$ and $p$ a prime. Suppose that for each $i\geq a$, $G$
satisfies~{\rm (\ref{commsthreetriv})}. 
If $N\geq (c-2){\rm ord}_p(c!)+a$, then
$$[z^{p^N},y]=[z,y]^{p^N}=[z,y^{p^N}].$$}

\proof By \ref{expofthirdterm}, $\langle z,y\rangle_3$ is of exponent
$p^{(c-3){\rm ord}_p(c!)+a}$. Let 
$$N\geq (c-2){\rm ord}_p(c!)+a.$$
Apply \ref{struiklemhtwo} to $[z^{p^N},y]$ and $[z,y^{p^N}]$, and use the
fact that $\langle z,y\rangle_3$ is of exponent $p^{(c-3){\rm
ord}_p(c!)+a}$ to obtain the~result.\endproof

\rmrk{strongehyp} Note that (\ref{commsthreetriv}) is a stronger
hypothesis than~(\ref{condonzandywithp}); but that if we set $p>c$ and
$a=1$, \ref{transferifenough} yields the same conclusion as
\ref{fromptoengelcondgeneral}.

\thm{generalniltwoinc}{Let $c>1$, $p$ a prime. Let $G\in {\cal N}_2$,
and let $H$ be a subgroup of~$G$.  Let $a\geq 0$ and let 
$N\geq (c-2){\rm ord}_p(c!) + a$.
If $x^{p^N}, y^{p^N}\in H$, then $[x,y]^{p^{2N-a}}\in {\rm
dom}_G^{{\cal N}_c}(H).$}

\proof Let $F= G \amalg_H^{{\cal N}_c} G$, and let $(\lambda,\rho)$ be
the universal pair of maps of~$G$ into~$F$. For simplicity, write
$\lambda(x)=x$, $\lambda(y)=y$, $\rho(x)=z$ and $\rho(y)=w$. 

By hypothesis, $x^{p^N}=z^{p^N}$ and $y^{p^N}=w^{p^N}$.
Note that if $i\geq N$, then
$$[z,y^{p^i},z] = [z,(y^{p^N})^{p^{i-N}},z]
= [z,(w^{p^N})^{p^{i-N}},z]
= [z,w^{p^i},z]
= e,$$
since $\rho(G)$ is a nilpotent group of class at
most~two. Analogously, for every $i\geq N$, we have
$[z^{p^i},y,y]=e$. 

We also have that $$[x,y]^{p^{2N-a}} = [x^{p^{2N-a}},y] = [z^{p^{2N-a}},y]$$
and $$[z,w]^{p^{2N-a}} = [z,w^{p^{2N-a}}] = [z,y^{p^{2N-a}}],$$
so it suffices to show that $[z^{p^{2N-a}},y]=[z,y^{p^{2N-a}}]$ in~$F$.

Also, we have that
$$2N-a = N + N-a
\geq ((c-2){\rm ord}_p(c!)+a) +N -a 
= (c-2){\rm ord}_p(c!)+N.$$

By \ref{transferifenough}, it follows that
$$[z^{p^{2N-a}},y]=[z,y]^{p^{2N-a}} = [z,y^{p^{2N-a}}]$$
in~$F$; therefore $[x,y]^{p^{2N-a}}\in {\rm dom}_G^{{\cal
N}_c}(H)$, as~claimed.\endproof

To finish this section we present an example, suggested by George
Bergman, of a nontrivial
dominion of a subgroup of an ${\cal N}_2$-group in~${\cal N}\!il$.

{\bf Example \numbeq{exampleofnontriv}.} Let $\Z[{1\over p}]$ be the
subring of~$\Q$ generated by $\Z$ and ${1\over p}$. Let $G$ be the
group with underlying set
$$\Z[{\textstyle{1\over p}}] \oplus \Z[{\textstyle {1\over p}}] \oplus
\Z[{\textstyle{1\over p}}]$$
and multiplication
$(a,b,c)\cdot (x,y,z) = (a+x, b+y, c+z+bx)$. This group is
isomorphic to the multiplicative group of upper triangular special
$3\times 3$ matrices over $\Z[{1\over p}]$, via the identification
$$(x,y,z) \longleftrightarrow \pmatrix{1&y&z\cr
                                       0&1&x\cr
                                       0&0&1\cr}$$

Let $H$ be the subgroup generated by $(1,0,0)$ and~$(0,1,0)$.  It is
not hard to verify that $H$ is the subgroup with underlying set
$\Z\times\Z\times\Z$, and~that ${\rm dom}_G^{{\cal N}_2}(H)$ is the
subgroup with underlying set $\Z\times \Z\times
\Z[{\textstyle{1\over p}}]$.

We claim that the dominion of~$H$ in~$G$ in the category ${\cal
N}\!il$ is also equal to the subgroup with underlying set
$\Z\times\Z\times \Z[{1\over p}]$.  For this, it suffices to show that
for all $c>1$ and all~$i>0$,
$$\left(0,0,{\textstyle{1\over p^i}}\right)\in {\rm dom}_G^{{\cal N}_c}(H).$$

Fix $c>1$ and $i>0$. Let $N = (c-2){\rm ord}_p(c!) + i$ and let
$x=({1\over p^N},0,0)$, $y=(0,-{1\over p^N},0)$.

Then $x^{p^N}=(1,0,0)\in H$, and $y^{p^N}=(0,-1,0)\in H$. By
\ref{generalniltwoinc}, 
$$[x,y]^{p^{2N-i}}\in {\rm dom}_G^{{\cal N}_c}(H).$$

Calculating directly we have that $[x,y]=(0,0,{1\over p^{2N}})$, and
hence
$$[x,y]^{p^{2N-i}}=\left(0,0,{\textstyle{1\over p^{i}}}\right)\in {\rm
dom}_G^{{\cal N}_c}(H),$$ as~claimed. 

Therefore, $H\propcont {\rm dom}_G^{{\cal N}_2}(H)={\rm dom}_G^{{\cal
N}_c}(H)$ for all $c>1$, so
$$H\propcont {\rm dom}_G^{{\cal N}_2}(H)={\rm dom}_G^{{\cal
N}\!il}(H).$$
In particular, there are instances of nontrivial dominions in~${\cal
N}\!il$.\endproof

Note that in Example~\ref{exampleofnontriv}, $G$ lies in ${\cal N}_2$;
since the dominions of~$H$ in~${\cal
N}_2$ and in~${\cal N}\!il$ are equal and properly contain~$H$, it
follows that the dominion of~$H$ in~$G$ in the category ${\cal
A}^2\cap {\cal N}\!il$ is also~nontrivial. Hence, there are nontrivial
dominions of subgroups of ${\cal N}_2$-groups in the category ${\cal
A}^2\cap {\cal N}\!il$. 

Also worthy of note is the fact that even though $H$ itself is finitely
generated, ${\rm dom}_G^{{\cal N}_2}(H)$ is not finitely
generated. Note as well that~$G$ is not finitely generated.  In fact,
one can prove that dominions of subgroups of finitely generated
nilpotent groups (of any class) are trivial in ${\cal N}\!il$, and
that dominions of subgroups of finitely generated ${\cal N}_2$-groups
are trivial in ${\cal A}^2\cap {\cal N}\!il$.
For the proofs of these assertions, see~{\bf [\cite{fgnilprelim}]}. 

\Section{Dominions in other varieties of nilpotent groups}{domsinothernil}

In this section, we will mention briefly some results which give
information on dominions
in other varieties of nilpotent groups.

First, it is known that any nonabelian torsion
free locally nilpotent variety~${\cal V}$ (that is, the relatively
free groups in~${\cal V}$ are torsion free, and every finitely
generated group in~${\cal V}$ is nilpotent) contains~${\cal N}_2$. See
for example~{\bf [\cite{thirtynine}]}. Therefore, we have:

\thm{domlocniltorsionfree}{Every nonabelian torsion free
locally nilpotent variety ${\cal V}$ has instances of
nontrivial~dominions. Namely, if $p$ is a prime greater than the
nilpotency class of the relatively free ${\cal V}$ group on four
generators, then $$
\langle x^p,y^p\rangle\propcont{\rm dom}_{F^{(2)}(x,y)}^{\cal V}\bigl(\langle
x^p,y^p\rangle) =
\langle x^p,y^p,[x,y]^p\rangle.$$}
 
\proof Let $c$ be the nilpotency class of the relatively free ${\cal
V}$ group on four generators. As we noted above, ${\cal V}$
will contain $F^{\scriptstyle(2)}(x,y)$.

For every subgroup~$H$ of~$F^{\scriptscriptstyle (2)}(x,y)$, the
amalgamated~coproduct $$F^{\scriptscriptstyle (2)}(x,y)\amalg_H^{\cal
V} F^{\scriptscriptstyle (2)}(x,y)$$ is generated by four elements,
and hence is a quotient of the relatively free ${\cal V}$-group of
rank four. In particular, it is nilpotent of class at most~$c$.

Let $p$ be a prime, $p\geq c$. If ${\cal V'}$ is the (nilpotent)
variety generated by the relatively free ${\cal V}$ group on four
generators, then
$$F^{\scriptstyle (2)}(x,y)\amalg_H^{\cal V} F^{\scriptstyle (2)}(x,y)
= F^{\scriptstyle (2)}(x,y)\amalg_H^{\cal V'} F^{\scriptstyle
(2)}(x,y);$$
hence $${\rm dom}_{F^{\scriptstyle (2)}(x,y)}^{\cal V}\bigl(\langle
x^p,y^p\rangle\bigr) = {\rm dom}_{F^{\scriptstyle (2)}(x,y)}^{\cal
V'}\bigl(\langle x^p,y^p\rangle\bigr)$$
which by \ref{lotsofvarsbetween} equals $\langle
x^p,y^p,[x,y]^p\rangle$, as~claimed.\endproof

In the case of the varieties ${\cal B}_p\cap{\cal N}_c$,
that is, nilpotent groups of class at most~$c$ and exponent~$p$, 
with $p$ a
prime greater than $c$, (the latter condition is sometimes called
{\it of small class}), the answer is provided by a theorem of Maier
{\bf [\cite{nilexpp}]}, which implies that in these varieties
dominions are~trivial. Specifically, given two groups $G$ and~$K$ in
the variety in question, and an amalgam $[G,K;H]$,
Maier obtains necessary and
sufficient conditions for the weak embeddability of the amalgam into a
group~$M$ in the~variety. He then~proves:
 
\thm{amalgforp}{{\rm (Maier, Corollary~1.3 in~{\bf [\cite{nilexpp}]})}
Suppose that $c<p$, and let $G,K\in {\cal B}_p\cap {\cal N}_c$. If
$[G,K;H]$ is an amalgam, and is weakly embeddable in a group $M\in
{\cal B}_p\cap {\cal N}_c$, then it is strongly embeddable in a group
$N$ in ${\cal B}_p\cap {\cal N}_c$.\noproof}
 
Since a special amalgam $[G,G;H]$ is always weakly embeddable, we
conclude~that:
 
\cor{specialamalgforp}{Suppose that $c<p$. Let $G\in{\cal B}_p\cap
{\cal N}_c$, and let $H$ be a subgroup of~$G$. Then ${\rm
dom}_{G}^{{\cal B}_p\cap {\cal N}_c}({H})=H$.\noproof}

Finally, we look at the nonabelian varieties of nilpotent groups of
class at most two.
First we state the classification of subvarieties of~${\cal N}_2$. It
follows as a corollary from the classification of subvarieties
of~${\cal N}_3$, due to J\'onsson and Remeslennikov (see~{\bf
[\cite{jonsson}]} and~{\bf [\cite{classifthree}]}).

\thm{classtwo}{Every variety of nilpotent groups of class at most~2
may be defined by the identities
$$x^m = [x_1,x_2]^n = [x_1,x_2,x_3] = e$$
for unique nonnegative integers $m$ and~$n$ satisfying
\hbox{$n|m/{\rm gcd}(2,m)$}, yielding a bijection between pairs of
nonnegative integers $(m,n)$ satisfying this condition, and varieties
of nilpotent groups of class at most~two.\noproof}

\thm{nnotsquarefree}{Let ${\cal V}$ be a variety of nilpotent groups
of class~$2$ corresponding to the pair $(m,n)$, and suppose that $n$
is not square free. Then there are nontrivial dominions in~${\cal
V}$.}

\proof If $m=n=0$, then ${\cal V}={\cal N}_2$ and there is nothing to
prove. Otherwise, let~$G$ be the relatively free group of rank~2
in~${\cal V}$, generated by~$x$ and~$y$. Let $p$ be a prime number
such that $p^2|n$. Let $H=\langle x^p,y^p\rangle$. If $m>0$, it is not
hard to verify that
$$H = \Bigl\{ x^{ap}y^{bp}[x,y]^{cp^2}\,\Bigm|\, 0\leq pa,pb< m,\quad
0\leq cp^2<n\Bigr\}.$$
However, $[x,y]^p\in {\rm dom}_G^{\cal V}(H)$, so $H\propcont {\rm dom}^{\cal
V}_G(H)$, as~claimed. If $m=0$, we proceed as above noting that $a$
and~$b$ may now be any integers, instead of being bounded by~$m/p$. This
proves the~theorem.\endproof

What about other varieties of nilpotent groups of class~two? We can
settle the matter for a few of the remaining varieties.

For $m>0$, we may use the facts that a finite nilpotent group is the
direct product of its Sylow subgroups, that a finitely generated
torsion nilpotent group is finite, and that dominions respect finite
direct products, to reduce to the case where $m$ is a prime power.
If $m=p$ a prime, then we are in the situation of
\ref{specialamalgforp}, so dominions are trivial. This leaves the
varieties $(p^n,p)$ with $n>1$, and $(0,n)$ with $n$ square free still
open. We ask:

{\bf Question \numbeq{squarefreenquestion}.} Are dominions trivial in
the subvarieties of ${\cal N}_2$ corresponding to pairs of positive
integers $(p^a,p)$ with $p$ a prime, $a>1$? Are dominions trivial in
the subvarieties corresponding to pairs $(0,n)$ with $n>1$ square~free?

My guess is that dominions will be trivial in those varieties, but I
am at present unable to prove this~guess.

\vskip\baselineskip
\centerline{AKNOWLEDGEMENTS}

The author was supported in part by a fellowship from the Programa de
Formaci\'on y Superaci\'on del Personal Acad\'emico de la UNAM,
administered by the DGAPA.

{\baselineskip=10pt
\smallheadfont
\def\bf{\smallheadbf}
\def\it{\smallheadit}
%
\ifnum0<\citations{\par\bigbreak
\filbreak\centerline{\rm REFERENCES}\par\frenchspacing}\fi
%
\ifundefined{xthreeNB}\else
\item{\bf [\refer{threeNB}]}{G{.} Baumslag, B.H.~Neumann,
H.~Neumann, and P.M.~Neumann. {On varieties generated by a
finitely generated group,} {\it Math.\ Z.} {\bf 86} (1964)
\hbox{93--122}. {MR:30\#138}}\par\filbreak\fi
\ifundefined{xbergman}\else
\item{\bf [\refer{bergman}]}{George M.~Bergman. {``An Invitation to
General Algebra and Universal Constructions,''} {Berkeley
Mathematics Lecture Notes 7,} 1995.}\par\filbreak\fi
\ifundefined{xordersberg}\else
\item{\bf [\refer{ordersberg}]}{George M.~Bergman, {Ordering
coproducts of groups and semigroups,} {\it J.\ Algebra} {\bf 133} (1990)
no. 2, \hbox{313--339}. {MR:91j:06035}}\par\filbreak\fi
\ifundefined{xbirkhoff}\else
\item{\bf [\refer{birkhoff}]}{Garrett Birkhoff, {On the structure
of abstract algebras.} {\it Proc.\ Cambridge\ Philos.\ Soc.} {\bf
31} (1935), \hbox{433--454}.}\par\filbreak\fi
\ifundefined{xbrown}\else
\item{\bf [\refer{brown}]}{Kenneth S.~Brown, {``Cohomology of
Groups'' 2nd Edition,} {GTM 87\/},
Springer Verlag,~1994. {MR:96a:20072}}\par\filbreak\fi
\ifundefined{xmetab}\else
\item{\bf [\refer{metab}]}{O.N.~Golovin, {Metabelian products of
groups,}
{\it Amer.\ Math.\ Soc.\ Translations series 2}, {\bf 2} (1956),
\hbox{117--131.} {MR:17,824b}}\par\filbreak\fi
\ifundefined{xhall}\else
\item{\bf [\refer{hall}]}{M.~Hall, {``The Theory of Groups''}
Mac~Millan Company,~1959. {MR:21\#1996}}\par\filbreak\fi
\ifundefined{xphall}\else
\item{\bf [\refer{phall}]}{P.~Hall, {Verbal and marginal
subgroups,} {\it J.\ Reine\ Angew.\ Math.\/} {\bf 182} (1940)
\hbox{156--157.} {MR:2,125i}}\par\filbreak\fi
\ifundefined{xheineken}\else
\item{\bf [\refer{heineken}]}{H.~Heineken, {Engelsche Elemente der
L\"ange drei,} {\it Illinois J.\ of Math.} {\bf 5} (1961)
\hbox{681--707.} {MR:24\#A1319}}\par\filbreak\fi
\ifundefined{xherman}\else
\item{\bf [\refer{herman}]}{Krzysztof~Herman, {Some remarks on
the twelfth problem of Hanna Neumann,} {\it Publ.\ Math.\ Debrecen}
{\bf 37} (1990)  no. 1--2, \hbox{25--31.} {MR:91f:20030}}\par\filbreak\fi
\ifundefined{xherstein}\else
\item{\bf [\refer{herstein}]}{I.~N. Herstein, {``Topics in
Algebra,''} Blaisdell Publishing Co.,~1964.}\par\filbreak\fi
\ifundefined{xepisandamalgs}\else
\item{\bf [\refer{episandamalgs}]}{Peter M.~Higgins, {Epimorphisms
and amalgams,} {\it
Colloq.\ Math.} {\bf 56} no.~1 (1988) \hbox{1--17.}
{MR:89m:20083}}\par\filbreak\fi
\ifundefined{xhigmanpgroups}\else
\item{\bf [\refer{higmanpgroups}]}{Graham Higman, {Amalgams of
$p$-groups,} {\it J. of~Algebra} {\bf 1} (1964)
\hbox{301--305.} {MR:29\#4799}}\par\filbreak\fi
\ifundefined{xhigmanremarks}\else
\item{\bf [\refer{higmanremarks}]}{Graham Higman, {Some remarks
on varieties of groups,} {\it Quart.\ J.\ of Math.\ (Oxford) (2)} {\bf
10} (1959), \hbox{165--178.} {MR:22\#4756}}\par\filbreak\fi
\ifundefined{xhughes}\else
\item{\bf [\refer{hughes}]}{N.J.S.~Hughes, {The use of bilinear
mappings in the classification of groups of class~$2$,} {\it Proc.\
Amer.\ Math.\ Soc.\ } {\bf 2} (1951) \hbox{742--747.}
{MR:13,528e}}\par\filbreak\fi
\ifundefined{xisbelltwo}\else
\item{\bf [\refer{isbelltwo}]}{J.~M.~Howie and J.~R. Isbell, {
Epimorphisms and dominions II,} {\it J.\ Algebra {\bf
6}}(1967) \hbox{7--21.} {MR:35\#105b}}\par\filbreak\fi
\ifundefined{xisbellone}\else
\item{\bf [\refer{isbellone}]}{J. R. Isbell, {Epimorphisms and
dominions} {\it in} { 
``Proc.~of the Conference on Categorical Algebra, La Jolla 1965,''\/}
Lange and Springer, New
York~1966. MR:35\#105a (The statement of the
Zigzag Lemma for {\it rings} in this paper is incorrect. The correct
version is stated in~{\bf [\cite{isbellfour}]}.)}\par\filbreak\fi
\ifundefined{xisbellthree}\else
\item{\bf [\refer{isbellthree}]}{J. R. Isbell, {Epimorphisms and
dominions III,} {\it Amer.\ J.\ Math.\ }{\bf 90} (1968)
\hbox{1025--1030.} {MR:38\#5877}}\par\filbreak\fi
\ifundefined{xisbellfour}\else
\item{\bf [\refer{isbellfour}]}{J. R. Isbell, {Epimorphisms and
dominions IV,} {\it J.\ London Math.\ Soc.~(2),}
{\bf 1} (1969) \hbox{265--273.} {MR:41\#1774}}\par\filbreak\fi
\ifundefined{xjones}\else
\item{\bf [\refer{jones}]}{Gareth A.~Jones, {Varieties and simple
groups,} {\it J.\ Austral.\ Math.\ Soc.} {\bf 17} (1974)
\hbox{163--173.} {MR:49\#9081}}\par\filbreak\fi
\ifundefined{xjonsson}\else
\item{\bf [\refer{jonsson}]}{B.~J\'onsson, {Varieties of groups of
nilpotency three,} {\it Notices Amer.\ Math.\ Soc.} {\bf 13} (1966)
488.}\par\filbreak\fi
\ifundefined{xwreathext}\else
\item{\bf [\refer{wreathext}]}{L.~Kaloujnine and Marc Krasner,
{Produit complet des groupes de permutations et le probl\`eme
d'extension des groupes III,} {\it Acta Sci.\ Math.\ Szeged} {\bf 14}
(1951) \hbox{69--82}. {MR:14,242d}}\par\filbreak\fi
\ifundefined{xkhukhro}\else
\item{\bf [\refer{khukhro}]}{Evgenii I. Khukhro, {``Nilpotent Groups
and their Automorphisms,''} {de Gruyter Expositions in Mathematics}
{\bf 8}, New York 1993. {MR:94g:20046}}\par\filbreak\fi
\ifundefined{xkleimanbig}\else
\item{\bf [\refer{kleimanbig}]}{Yu.~G. Kle\u{\i}man, {On
identities in groups,} {\it Trans.\ Moscow Math.\ Soc.\ } 1983,
Issue 2, \hbox{63--110}. {MR:84e:20040}}\par\filbreak\fi
\ifundefined{xthirtynine}\else
\item{\bf [\refer{thirtynine}]}{L. G. Kov\'acs, {The thirty-nine
varieties,} {\it Math.\ Scientist} {\bf 4} (1979)
\hbox{113--128.} {MR:81m:20037}}\par\filbreak\fi
\ifundefined{xlamssix}\else
\item{\bf [\refer{lamssix}]}{T.Y. Lam and David B. Leep, {
Combinatorial structure on the automorphism group of~$S_6$,} {\it
Expo. Math.} {\bf 11} (1993) \hbox{289--308.}
{MR:94i:20006}}\par\filbreak\fi
\ifundefined{xlevione}\else
\item{\bf [\refer{levione}]}{F.~W. Levi, {Groups on which the
commutator relation 
satisfies certain algebraic conditions,} {\it J.\ Indian Math.\ Soc.\ New
Series} {\bf 6}(1942), \hbox{87--97.} {MR:4,133i}}\par\filbreak\fi
\ifundefined{xgermanlevi}\else
\item{\bf [\refer{germanlevi}]}{F.~W. Levi and B. L. van der Waerden,
{\"Uber eine 
besondere Klasse von Gruppen,} {\it Abhandl.\ Math.\ Sem.\ Univ.\ Hamburg}
{\bf 9}(1932), \hbox{154--158.}}\par\filbreak\fi
\ifundefined{xlichtman}\else
\item{\bf [\refer{lichtman}]}{A. L. Lichtman, {Necessary and
sufficient conditions for the residual nilpotence of free products of
groups,} {\it J. Pure and Applied Algebra} {\bf 12} no. 1 (1978),
\hbox{49--64.} {MR:58\#5938}}\par\filbreak\fi
\ifundefined{xmaxofan}\else
\item{\bf [\refer{maxofan}]}{Martin W. Liebeck, Cheryl E. Praeger, 
and Jan Saxl, {A classification of the maximal subgroups of the
finite alternating and symmetric groups,} {\it J. of Algebra} {\bf
111}(1987), \hbox{365--383.} {MR:89b:20008}}\par\filbreak\fi
\ifundefined{xepisingroups}\else
\item{\bf [\refer{episingroups}]}{C.E. Linderholm, {A group
epimorphism is surjective,} {\it Amer.\ Math.\ Monthly\ }77
\hbox{176--177.}}\par\filbreak\fi
\ifundefined{xmckay}\else
\item{\bf [\refer{mckay}]}{Susan McKay, {Surjective epimorphisms
in classes
of groups,} {\it Quart.\ J.\ Math.\ Oxford (2),\/} {\bf 20} (1969),
\hbox{87--90.} {MR:39\#1558}}\par\filbreak\fi
\ifundefined{xmachenry}\else
\item{\bf [\refer{machenry}]}{T. MacHenry, {The tensor product and
the 2nd nilpotent product of groups,} {\it Math. Z.\/} {\bf 73}
(1960), \hbox{134--145.} {MR:22\#11027a}}\par\filbreak\fi

\ifundefined{xmaclane}\else
\item{\bf [\refer{maclane}]}{Saunders Mac Lane, {``Categories for
the Working Mathematician,''} {GTM 5},
Springer Verlag (1971). {MR:50\#7275}}\par\filbreak\fi
\ifundefined{xbilinearprelim}\else
\item{\bf [\refer{bilinearprelim}]}{Arturo Magidin, {Bilinear maps
and central extensions of abelian groups,} {\it Submitted.}}\par\filbreak\fi
\ifundefined{xprodvarprelim}\else
\item{\bf [\refer{prodvarprelim}]}{Arturo Magidin, {Dominions in decomposable
varieties of groups,} {\it Submitted.}}\par\filbreak\fi
\ifundefined{xmythesis}\else
\item{\bf [\refer{mythesis}]}{Arturo Magidin, {``Dominions in
Varieties of Groups,''} Doctoral dissertation, University of
California at Berkeley, May 1998.}\par\filbreak\fi
\ifundefined{xnildomsprelim}\else
\item{\bf [\refer{nildomsprelim}]}{Arturo Magidin, {Dominions in
varieties of nilpotent groups,} {\it Comm.\ Alg.} to appear.}\par\filbreak\fi
\ifundefined{xsimpleprelim}\else
\item{\bf [\refer{simpleprelim}]}{Arturo Magidin, {Dominions in
varieties generated by simple groups,} {\it Submitted.}}\par\filbreak\fi
\ifundefined{xdomsmetabprelim}\else
\item{\bf [\refer{domsmetabprelim}]}{Arturo Magidin, {Dominions
in the variety of metabelian groups,}
{\it Submitted.}}\par\filbreak\fi
\ifundefined{xfgnilprelim}\else
\item{\bf [\refer{fgnilprelim}]}{Arturo Magidin, {Dominions of
finitely generated nilpotent groups,} {\it
Comm.\ Algebra,} to appear}\par\filbreak\fi 
\ifundefined{xwordsprelim}\else
\item{\bf [\refer{wordsprelim}]}{Arturo Magidin, {
Words and dominions,} {\it Submitted.}}\par\filbreak\fi
\ifundefined{xabsclosed}\else
\item{\bf [\refer{absclosed}]}{Arturo Magidin, {Absolutely closed
nil-2 groups,} {\it Submitted.}}\par\filbreak\fi
\ifundefined{xmagnus}\else
\item{\bf [\refer{magnus}]}{Wilhelm Magnus, Abraham Karras, and
Donald Solitar, {``Combinatorial Group Theory,''} 2nd Edition; Dover
Publications, Inc.~1976. {MR:53\#10423}}\par\filbreak\fi
\ifundefined{xamalgtwo}\else
\item{\bf [\refer{amalgtwo}]}{Berthold J. Maier, {Amalgame
nilpotenter Gruppen
der Klasse zwei II,} {\it Publ.\ Math.\ Debrecen} {\bf 33}(1986),
\hbox{43--52.} {MR:87k:20050}}\par\filbreak\fi
\ifundefined{xnilexpp}\else
\item{\bf [\refer{nilexpp}]}{Berthold J. Maier, {On nilpotent
groups of exponent $p$,} {\it J.\ Algebra} {\bf 127} (1989)
\hbox{279--289.} {MR:91b:20046}}\par\filbreak\fi
\ifundefined{xmaltsev}\else
\item{\bf [\refer{maltsev}]}{A. I. Maltsev, {Generalized
nilpotent algebras and their associated groups} (Russian), {\it
Mat.\ Sbornik N.S.} {\bf 25(67)} (1949) \hbox{347--366.} ({\it
Amer.\ Math.\ Soc.\ Translations Series 2} {\bf 69} 1968,
\hbox{1--21.}) {MR:11,323b}}\par\filbreak\fi
\ifundefined{xmaltsevtwo}\else
\item{\bf [\refer{maltsevtwo}]}{A. I. Maltsev, {Homomorphisms onto
finite groups} (Russian), {\it Ivanov. gosudarst. ped. Inst., u\v
cenye zap., fiz-mat. Nauk} {\bf 18} (1958)
\hbox{49--60.}}\par\filbreak\fi
\ifundefined{xmorandual}\else
\item{\bf [\refer{morandual}]}{S. Moran, {Duals of a verbal
subgroup,} {\it J.\ London Math.\ Soc.} {\bf 33} (1958)
\hbox{220--236.} {MR:20\#3909}}\par\filbreak\fi
\ifundefined{xhneumann}\else
\item{\bf [\refer{hneumann}]}{Hanna Neumann, {``Varieties of
Groups,''} {Ergebnisse der Mathematik und ihrer Grenz\-ge\-biete,\/}
New series, Vol.~37, Springer Verlag~1967. {MR:35\#6734}}\par\filbreak\fi
\ifundefined{xneumannwreath}\else
\item{\bf [\refer{neumannwreath}]}{Peter M.~Neumann, {On the
structure of standard wreath products of groups,} {\it Math.\
Zeitschr.\ }{\bf 84} (1964) \hbox{343--373.} {MR:32\#5719}}\par\filbreak\fi
\ifundefined{xpneumann}\else
\item{\bf [\refer{pneumann}]}{Peter M.~Neumann, {Splitting groups
and projectives
in varieties of groups,} {\it Quart.\ J.\ Math.\ Oxford} (2), {\bf
18} (1967),
\hbox{325--332.} {MR:36\#3859}}\par\filbreak\fi
\ifundefined{xoates}\else
\item{\bf [\refer{oates}]}{Sheila Oates, {Identical Relations in
Groups,} {\it J.\ London Math.\ Soc.} {\bf 38} (1963),
\hbox{71--78.} {MR:26\#5043}}\par\filbreak\fi
\ifundefined{xolsanskii}\else
\item{\bf [\refer{olsanskii}]}{A. Ju.~Ol'\v{s}anski\v{\i}, {On the
problem of a finite basis of identities in groups,} {\it
Izv.\ Akad.\ Nauk.\ SSSR} {\bf 4} (1970) no. 2
\hbox{381--389.}}\par\filbreak\fi
\ifundefined{xremak}\else
\item{\bf [\refer{remak}]}{R. Remak, {\"Uber minimale invariante
Untergruppen in der Theorie der end\-lichen Gruppen,} {\it
J.\ reine.\ angew.\ Math.} {\bf 162} (1930),
\hbox{1--16.}}\par\filbreak\fi
\ifundefined{xclassifthree}\else
\item{\bf [\refer{classifthree}]}{V. N. Remeslennikov, {Two
remarks on 3-step nilpotent groups} (Russian), {\it Algebra i Logika
Sem.} (1965) no.~2 \hbox{59--65.} {MR:31\#4838}}\par\filbreak\fi
\ifundefined{xrotman}\else
\item{\bf [\refer{rotman}]}{J.J. Rotman, {``Introduction to the Theory of
Groups 4th edition,''} {GTM~119},
Springer Verlag,~1994. {MR:95m:20001}}\par\filbreak\fi
\ifundefined{xsaracino}\else
\item{\bf [\refer{saracino}]}{D. Saracino, {Amalgamation bases for
nil-$2$ groups,} {\it Alg.\ Universalis\/} {\bf 16} (1983),
\hbox{47--62.} {MR:84i:20035}}\par\filbreak\fi
\ifundefined{xscott}\else
\item{\bf [\refer{scott}]}{W. R. Scott, {``Group Theory,''} Prentice
Hall,~1964. {MR:29\#4785}}\par\filbreak\fi
\ifundefined{xsmelkin}\else
\item{\bf [\refer{smelkin}]}{A. L. \v{S}mel'kin, {Wreath products and
varieties of groups} (Russian), {\it Dokl.\ Akad.\ Nauk S.S.S.R.\/} {\bf
157} (1964), \hbox{1063--1065} Transl.: {\it Soviet Math.\ Dokl.\ } {\bf
5} (1964), \hbox{1099--1011}. {MR:33\#1352}}\par\filbreak\fi
\ifundefined{xstruikone}\else
\item{\bf [\refer{struikone}]}{Ruth Rebekka Struik, {On nilpotent
products of cyclic groups,} {\it Canadian J.\ Math.}
{\bf 12} (1960)
\hbox{447--462}. {MR:22\#11028}}\par\filbreak\fi
\ifundefined{xstruiktwo}\else
\item{\bf [\refer{struiktwo}]}{Ruth Rebekka Struik, {On nilpotent
products of cyclic groups II,} {\it Canadian J.\
Math.\/} {\bf 13} (1961) \hbox{557--568.}
{MR:26\#2486}}\par\filbreak\fi
\ifundefined{xvlee}\else
\item{\bf [\refer{vlee}]}{M. R. Vaughan-Lee, {Uncountably many
varieties of groups,} {\it Bull.\ London Math.\ Soc.} {\bf 2} (1970)
\hbox{280--286.} {MR:43\#2054}}\par\filbreak\fi
\ifundefined{xweibel}\else
\item{\bf [\refer{weibel}]}{Charles Weibel, {``Introduction to
Homological Algebra,''} Cambridge University
Press~1994. {MR:95f:18001}}\par\filbreak\fi 
\ifundefined{xweigelone}\else
\item{\bf [\refer{weigelone}]}{T. S. Weigel, {Residual properties
of free groups,} {\it J.\ Algebra} {\bf 160} (1993)
\hbox{14--41.} {MR:94f:20058a}}\par\filbreak\fi
\ifundefined{xweigeltwo}\else
\item{\bf [\refer{weigeltwo}]}{T. S. Weigel, {Residual properties
of free groups II,} {\it Comm.\ in Algebra} {\bf 20}(5) (1992)
\hbox{1395--1425.} {MR:94f:20058b}}\par\filbreak\fi
\ifundefined{xweigelthree}\else 
\item{\bf [\refer{weigelthree}]}{T. S. Weigel, {Residual Properties
of free groups III,} {\it Israel J.\ Math.\ } {\bf 77} (1992)
\hbox{65--81.} {MR:94f:20058c}}\par\filbreak\fi
\ifundefined{xzstwo}\else
\item{\bf [\refer{zstwo}]}{Oscar Zariski and Pierre Samuel,
{``Commutative Algebra, Volume
II,''} Springer-Verlag~1976. {MR:52\#10706}}\par\filbreak\fi
\ifnum0<\citations\nonfrenchspacing\fi

}
\bigskip

\vfill
\eject
\immediate\closeout\aux
\end